\documentclass[12pt]{article}
\usepackage{amsmath,amsthm,amscd,amssymb}
\pdfoutput=1
\usepackage[margin=.75 in]{geometry}
\usepackage{subfig}
\usepackage{graphicx}
\usepackage[bookmarks=false]{hyperref}
\usepackage[labelfont=bf]{caption}
\allowdisplaybreaks
\usepackage{color}
\usepackage{cancel}
\usepackage{hyperref}
\theoremstyle{plain}

\theoremstyle{definition}

\theoremstyle{remark}

\newcommand{\mf}{\mathfrak}
\newcommand{\mc}{\mathcal}
\newcommand{\mb}{\mathbb}
\newcommand{\mbf}{\mathbf}

\newcommand{\Z}{\mathbb{Z}}

\newcommand{\R}{\mathbb{R}}
\newcommand{\C}{\mathbb{C}}
\newcommand{\ba}{\begin{aligned}}
\newcommand{\ea}{\end{aligned}}
\newcommand{\bt}{\begin{thm}}
\newcommand{\et}{\end{thm}}
\newcommand{\bc}{\begin{corollary}}
\newcommand{\ec}{\end{corollary}}
\newcommand{\bl}{\begin{lemma}}
\newcommand{\el}{\end{lemma}}
\newcommand{\bpf}{\begin{proof}}
\newcommand{\epf}{\end{proof}}
\newcommand{\bpb}{\begin{problem}}
\newcommand{\epb}{\end{problem}}
\newcommand{\bd}{\begin{definition}}
\newcommand{\ed}{\end{definition}}
\newcommand{\bn}{\begin{note}}
\newcommand{\en}{\end{note}}
\newcommand{\bp}{\begin{proposition}}
\newcommand{\ep}{\end{proposition}}
\newcommand{\be}{\begin{example}}
\newcommand{\ee}{\end{example}}
\newcommand{\bex}{\begin{exercise}}
\newcommand{\eex}{\end{exercise}}
\theoremstyle{plain}
\newtheorem{thm}{Theorem}[section]

\newtheorem{lemma}[thm]{Lemma}
\newtheorem{corollary}[thm]{Corollary}
\newtheorem{proposition}[thm]{Proposition}
\newtheorem{exercise}[thm]{Exercise}
\newtheorem{problem}[thm]{Problem}

\theoremstyle{definition}
\newtheorem{definition}[thm]{Definition}
\newtheorem{remark}[thm]{Remark}
\newtheorem{example}[thm]{Example}

\newtheorem*{convention}{Important Notational Convention}

\newtheorem {note}{Note}[section]
\swapnumbers
\theoremstyle{remark}

\theoremstyle{plain}

\title{Hard--edge asymptotics of the Jacobi growth process}
\author{ Mark Cerenzia and Jeffrey Kuan
\thanks{Electronic addresses: \texttt{cerenzia@princeton.edu}, \texttt{kuan@math.columbia.edu}} \\
}
\date{}
\begin{document}
\maketitle

\begin{abstract}

We introduce a two parameter ($\alpha, \beta>-1$) family of interacting particle systems with determinantal correlation kernels expressible in terms of Jacobi polynomials $\{ P^{(\alpha, \beta)}_k \}_{k \geq 0}$. The family includes previously discovered Plancherel measures for the infinite--dimensional orthogonal and symplectic groups. The construction uses certain BC--type orthogonal polynomials which generalize the characters of these groups.

The local asymptotics near the hard edge where one expects distinguishing behavior yields {\color{black} the multi--time $(\alpha, \beta)$--dependent discrete Jacobi kernel and the multi--time $\beta$--dependent hard--edge Pearcey kernel. For nonnegative integer values of $\beta$, the hard--edge Pearcey kernel had previously appeared in the asymptotics of non--intersecting squared Bessel paths at the hard edge.} 

\end{abstract}
In their study of the perturbed chiral GUE ensemble, Patrick Desrosiers and Peter J. Forrester \cite{pdpf1} (2008) introduce a process with a (fixed-time) kernel, depending on a parameter $\beta$, which can be written in terms of the Bessel function $J_{\beta}$ of the first kind. Later, in the asymptotic analysis of an interacting particle system near a hard edge with a reflecting barrier, Borodin--Kuan \cite{abjk1} (2009) find a multi--time kernel, which at fixed times corresponds to the $\beta=-1/2$ case of a $\beta$--dependent kernel arising in the work \cite{kf1} of Kuijlaars, Martinez-Finkelshtein, Wielonsky (2010) on non-intersecting squared Bessel paths at the hard edge at $0$. The authors conjectured that this kernel is the same as the one from \cite{pdpf1}. Indeed, Delvaux--Vet{\H{o}} \cite{sdbv1} (2014), {\color{black} assuming that $\beta$ is a nonnegative integer}, extends the results to a multi--time $\beta$--dependent kernel called the \textit{hard--edge Pearcey kernel}, and prove that it generalizes the kernel of \cite{pdpf1} for all $\beta>-1$ and the kernel of \cite{kf1} for nonnegative integer $\beta$. {\color{black} The terminology goes back to the Pearcey kernel, which had previously occurred at the soft edge of certain processes (see e.g. \cite{AMW},\cite{AOM},\cite{abk},\cite{abjk2},\cite{bh1},\cite{bh2},\cite{GZ},\cite{HHN},\cite{LW},\cite{aonr2},\cite{tw}).}

More recent work \cite{cerenzia1} constructs an interacting particle system similar to that of \cite{abjk1}, where jumps into the wall are suppressed instead of reflected. The resulting asymptotic kernel ends up being the $\beta=1/2$ case of the hard--edge Pearcey kernel. In light of the previous results, it is natural to ask for a $\beta$--dependent interacting particle system with asymptotics converging to the $\beta$--dependent hard--edge Pearcey kernel. This paper will construct a two parameter ($\alpha, \beta>-1$) family of interacting particle systems whose limiting behavior is the hard--edge Pearcey kernel. {\color{black} While the formula for the hard--edge Pearcey kernel had already been written in \cite{sdbv1}, the results here will establish that it is the correlation kernel of a well--defined determinantal point process for all $\beta>-1$, and not just for nonnegative integer $\beta$.}

 The methods in \cite{abjk1} and \cite{cerenzia1} use the representation theory of the orthogonal and symplectic groups, respectively. The connection to the parameter $\beta$ occurs through the characters, which can be written as a determinant of the Jacobi polynomials $P^{(\alpha,\beta)}$ for $(\alpha=\pm 1/2,\beta=-1/2)$ (in the orthogonal case) and $(\alpha=\pm 1/2,\beta=1/2)$ (in the symplectic case). Work of Okounkov and Olshanski \cite{aogo1} on BC--type orthogonal polynomials indicated that it should be possible to extend these {\color{black} methods} to general $\alpha,\beta>-1$, and it is through these Jacobi polynomials that the process is defined. See also \cite{def1}, \cite{def} and \cite{ww} for interacting particle systems with walls which are related to representation theory.

\begin{figure}[h!]
\centering
\hspace{5em}
\includegraphics[width=0.24\textwidth]{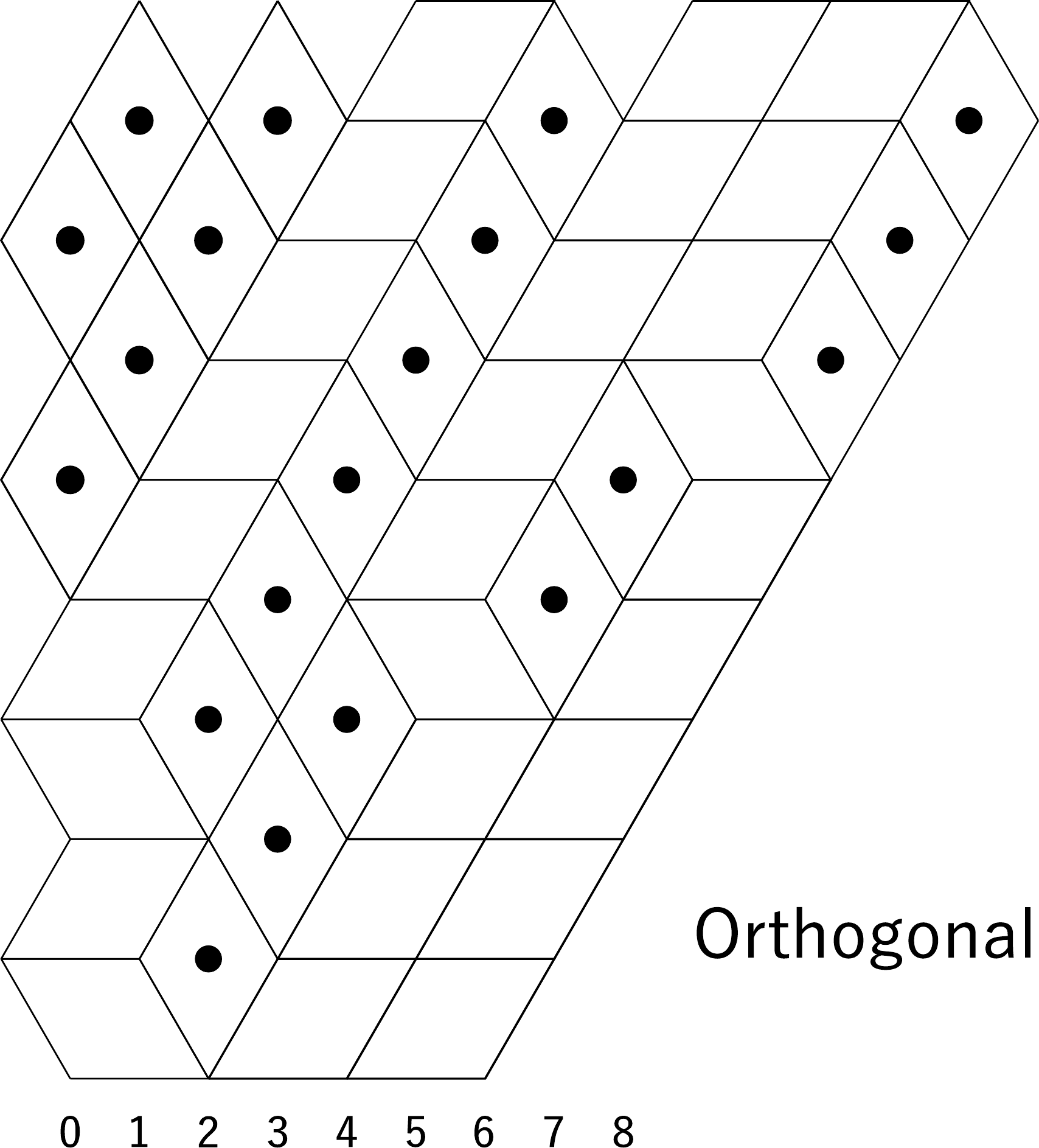}
\hspace{10em}
\includegraphics[width=0.24\textwidth]{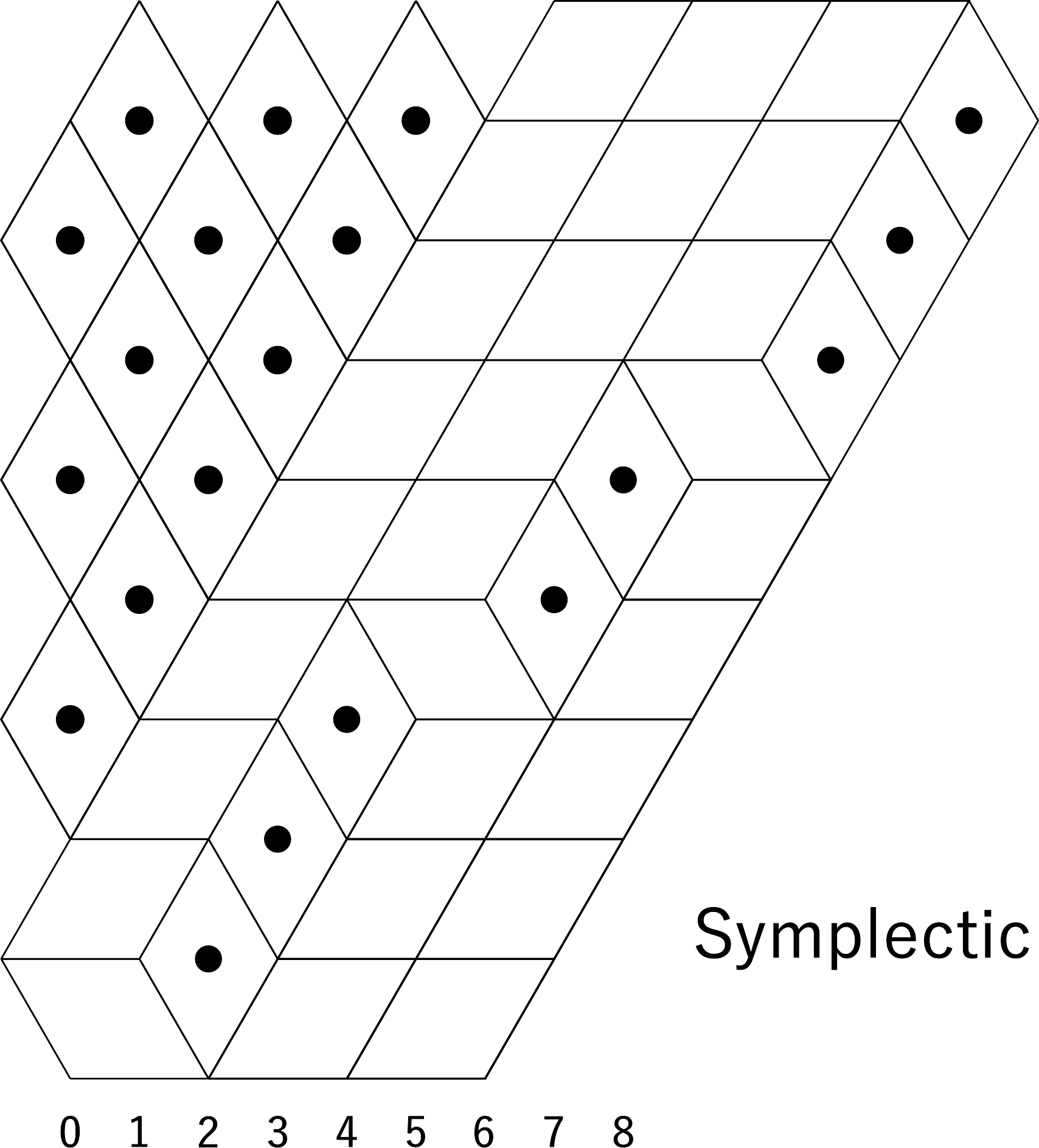}  
\hspace{5em}
\vspace{1.00mm}
\qquad
\newline
\centering
    \includegraphics[width=0.37\textwidth]{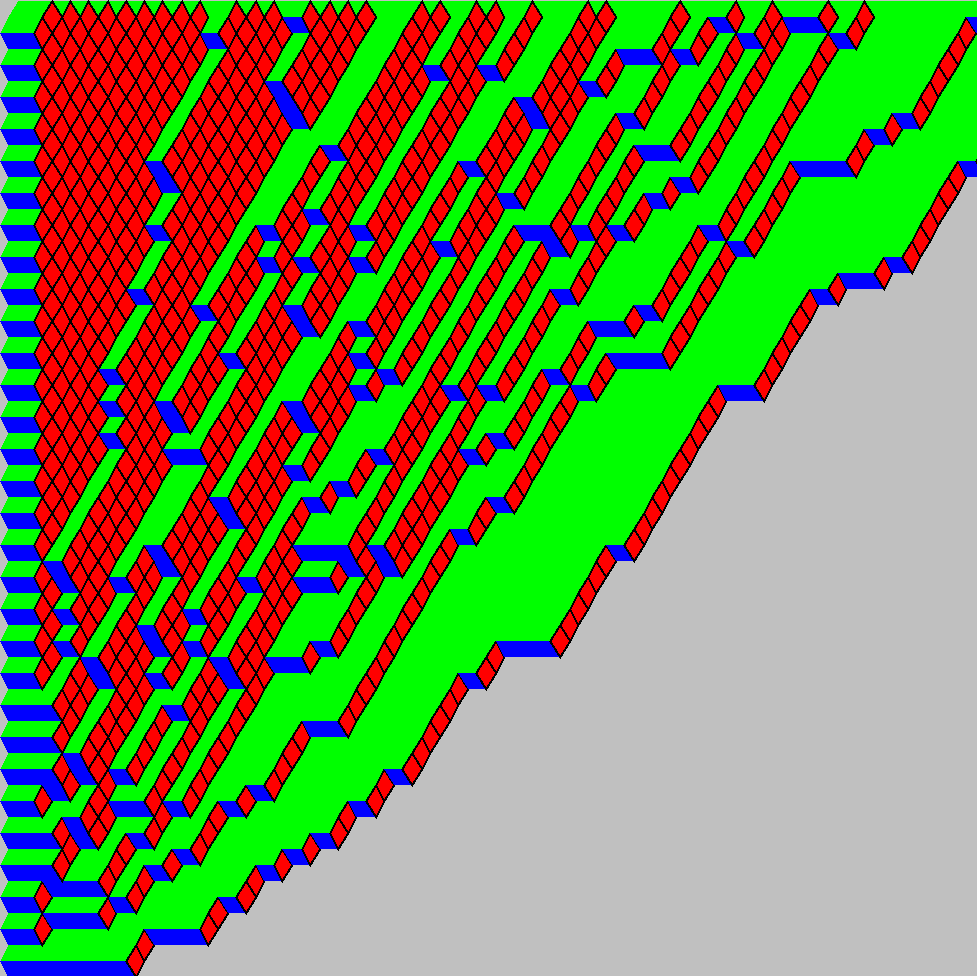}
    \hspace{5em}
    \includegraphics[width=0.37\textwidth]{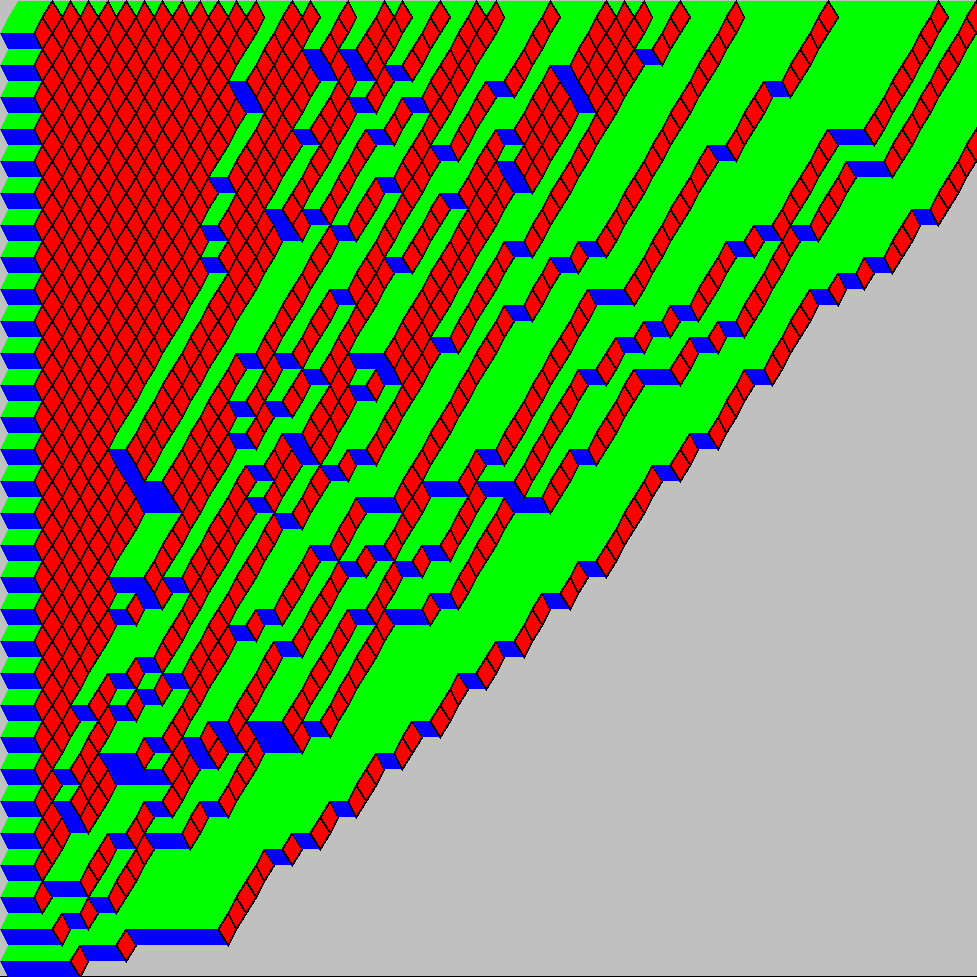}
\caption{The fourth step illustrates the distinctive wall behaviors: reflection (left) and suppression (right). Simulations of the orthogonal and {symplectic} cases indicate similar average behavior.} \label{modelcomparison}
\end{figure}

\begin{figure}[h] 
\centering
\includegraphics[width=0.4\textwidth]{symplecticlargetime.pdf} 
\hspace{.5em}
\includegraphics[width=0.327\textwidth]{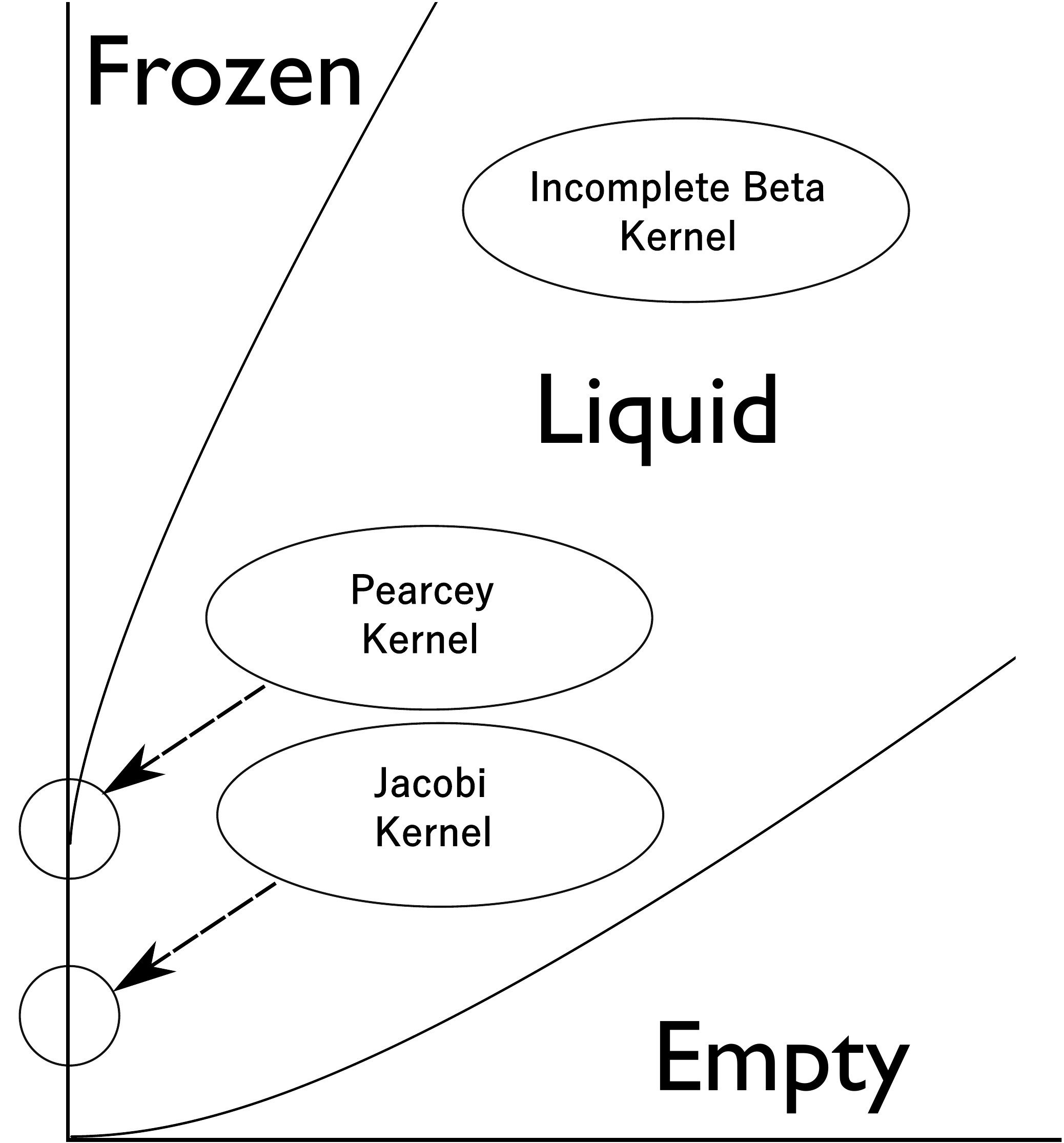} 
\caption{Summary of Borodin--Kuan \cite{abjk1} results. The hydrodynamic limit capturing average behavior indicates a densely-packed region of inactivity (``frozen") and a central region of activity (``liquid"), with the remaining area unreached (``empty"). Compare with Figure 4 of \protect\cite{abjk1}.} \label{averagebehavior}
\end{figure}

{\color{black} Let us now describe the model in more detail.} In orthogonal/{symplectic} Gelfand-Tsetlin patterns, level $n \in \Z_{> 0}:=\{1,2,3 \ldots \}$ has $r_n:=\lfloor (n+1)/2 \rfloor$ particles $\mc{X}^n_k(\gamma)$, $1 \leq k \leq r_n$, dynamically evolving according to the time parameter $\gamma \geq 0$ from the {\color{black} \textit{densely packed}} initial condition $\mc{X}^n_k(0) = r_n-k$. Each particle has two unit rate exponential clocks, one for rightward jumps and one for leftward jumps, all of which are independent. Suppose the right clock rings of $\mc{X}^n_k$. If $\mc{X}^n_k = \mc{X}^{n-1}_{k-1}$, then the right jump attempt is suppressed (a ``block" occurs). Otherwise, take the largest $c \geq 1$ such that
$$
\mc{X}^n_k = \mc{X}^{n+1}_k + (r_{n+1}-r_n) = \cdots = \mc{X}^{n+(c-1)}_k + \sum^{c-1}_{\ell=1} (r_{n+\ell+1}-r_{n+\ell}),
$$
then all $c$ of these particles jump to the right (a ``push" occurs). Similarly, leftward jump attempts by $\mc{X}^n_k$ can either be suppressed by a particle on the previous level $n-1$ in a state where it can push $\mc{X}^n_k$ by a right jump, or the left jump of $\mc{X}^n_k$ can push a particle on the next level $n+1$ whose rightward jumps would be blocked by $\mc{X}^n_k$. Lastly, we need to declare what occurs if $n \geq 1$ is odd and the left clock of a ``wall--particle" $\mc{X}^n_{r_n}$ rings when the particle is at position $0$ (i.e., at the wall). If such wall-jumps are reflected (i.e., a right jump is attempted), we call the dynamics for such point patterns \emph{orthogonal \ Plancherel \ Growth},  and \emph{{symplectic} \ Plancherel \ Growth} if such wall-jumps are suppressed. These dynamics can be depicted as follows (where the smaller arrows indicate which clock rings): 
\begin{equation} \label{maincoordinates}
 \ba
&| x^5_3 \ \ \ x^5_2 \ \  \  x^5_1 \\
&| x^4_2 \ \ \  x^4_1 \\
&| x^3_2 \ \ \  x^3_1\\
&|  x^2_1\to \\
&| x^1_1 
\ea
 \ \ \  \longrightarrow \ \ \ \ 
 \ba
&| x^5_3 \ \ \ x^5_2 \to \    x^5_1 \\
&|  x^4_2 \ \ \  \ \ \ \ x^4_1 \\
&| x^3_2 \ \ \  \ \ \  \ x^3_1\\
&|  \ \ \ \ \ x^2_1 \\
&|  x^1_1
\ea
 \ \ \  \longrightarrow \ \ \ \ 
 \ba
&| x^5_3 \ \ \  \ \ \  \ x^5_2 \  \   x^5_1 \\
&|  x^4_2 \ \  \leftarrow  x^4_1 \\
&| x^3_2 \ \ \  \ \ \ \ x^3_1\\
&|  \ \ \ \ \ x^2_1 \\
&| x^1_1
\ea
 \ \ \  \longrightarrow \ \ \ \ 
 \ba
&| x^5_3 \ \ x^5_2 \  \ \ \ \ \ \   x^5_1 \\
&|  x^4_2 \ \    x^4_1 \\
&| x^3_2 \ \ \  \ \ \ \ x^3_1\\
&|  \ \ \ \ \ x^2_1 \\
&| \overset{ \leftarrow }{x^1_1} 
\ea
 \ \ \  \longrightarrow \ \ \ \ \cdots 
\end{equation}
Note the blocking and pushing of a particle by straddling particles on the previous level serves to maintain \emph{interlacing conditions}, which in the chosen coordinates read as
\begin{equation}
\begin{cases}
x^{n+1}_{k+1} \leq  x^n_k < x^{n+1}_k, \ \ \ \text{if} \ \ n \ \ \text{even} \\
x^{n+1}_{k+1} \leq  x^n_k \leq x^{n+1}_k, \ \ \ \text{if} \ \ n \ \ \text{odd}
\end{cases}, \ \ \ 1\leq k \leq r_n.
\end{equation}
To make this visualization easier, see Figures \ref{modelcomparison} and \ref{averagebehavior}, which compares the dynamics in shifted coordinates. 
\begin{figure}[h!] 
\centering
\includegraphics[width=0.24\textwidth]{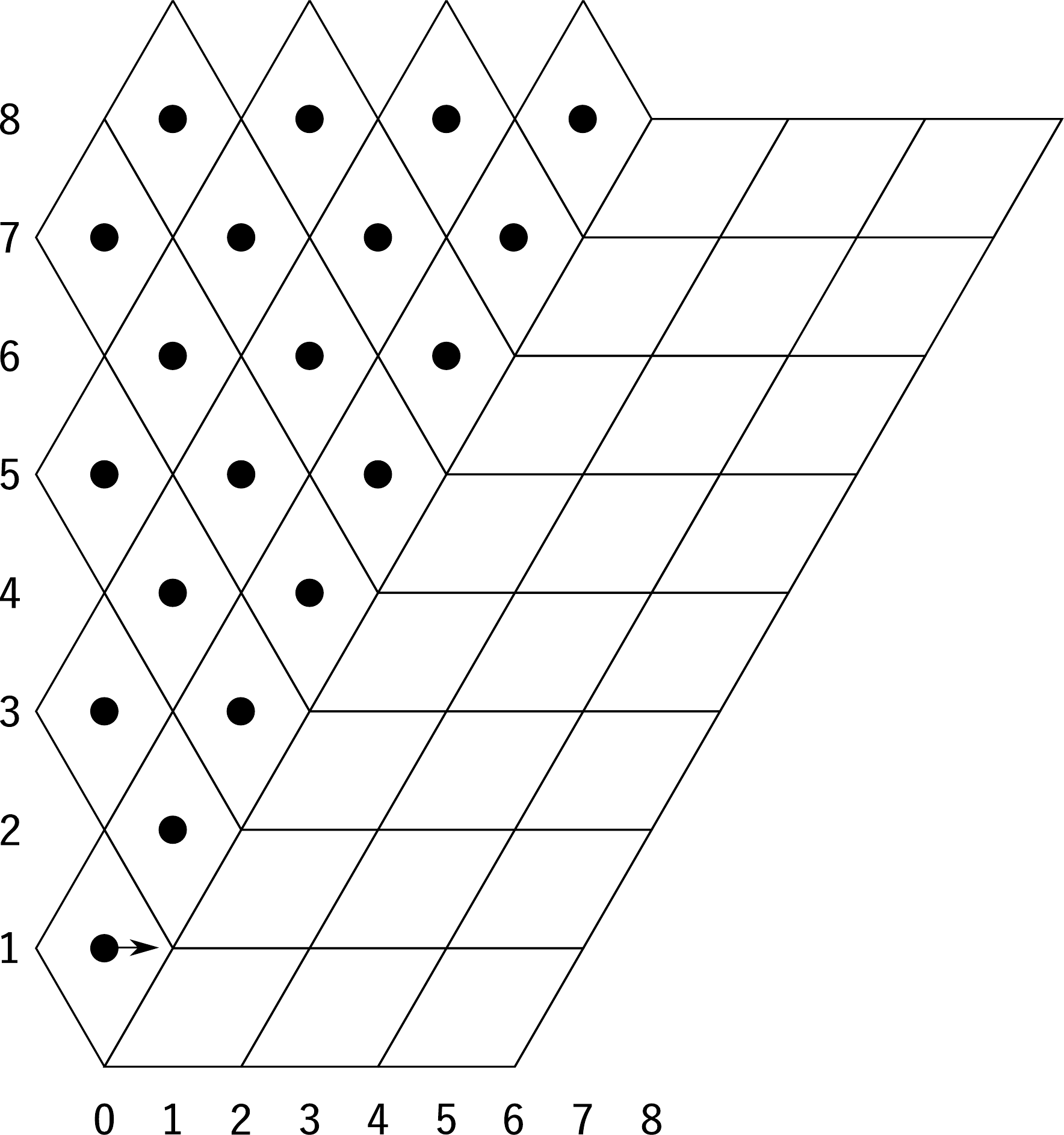} 
\includegraphics[width=0.24\textwidth]{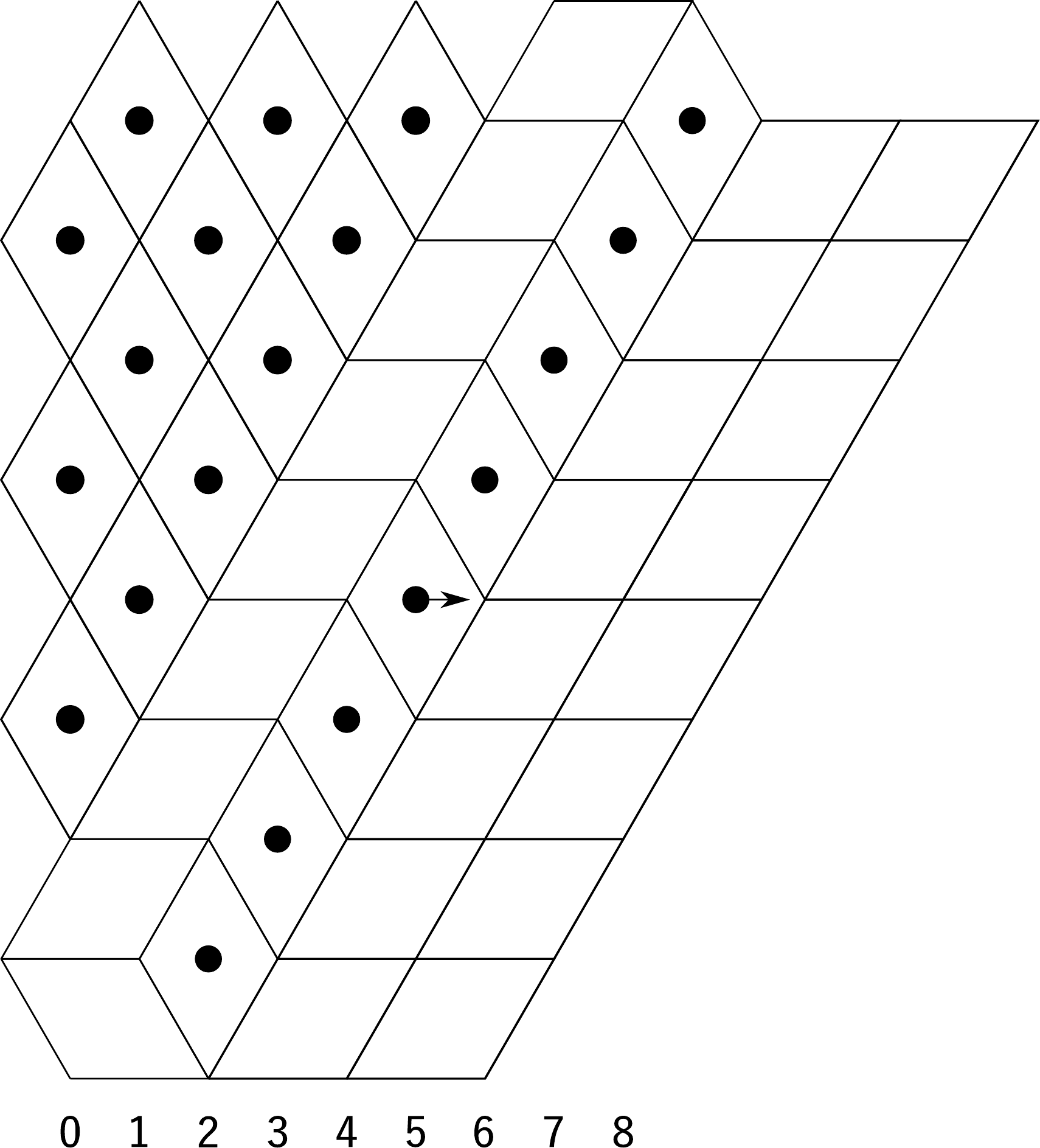} 
\includegraphics[width=0.24\textwidth]{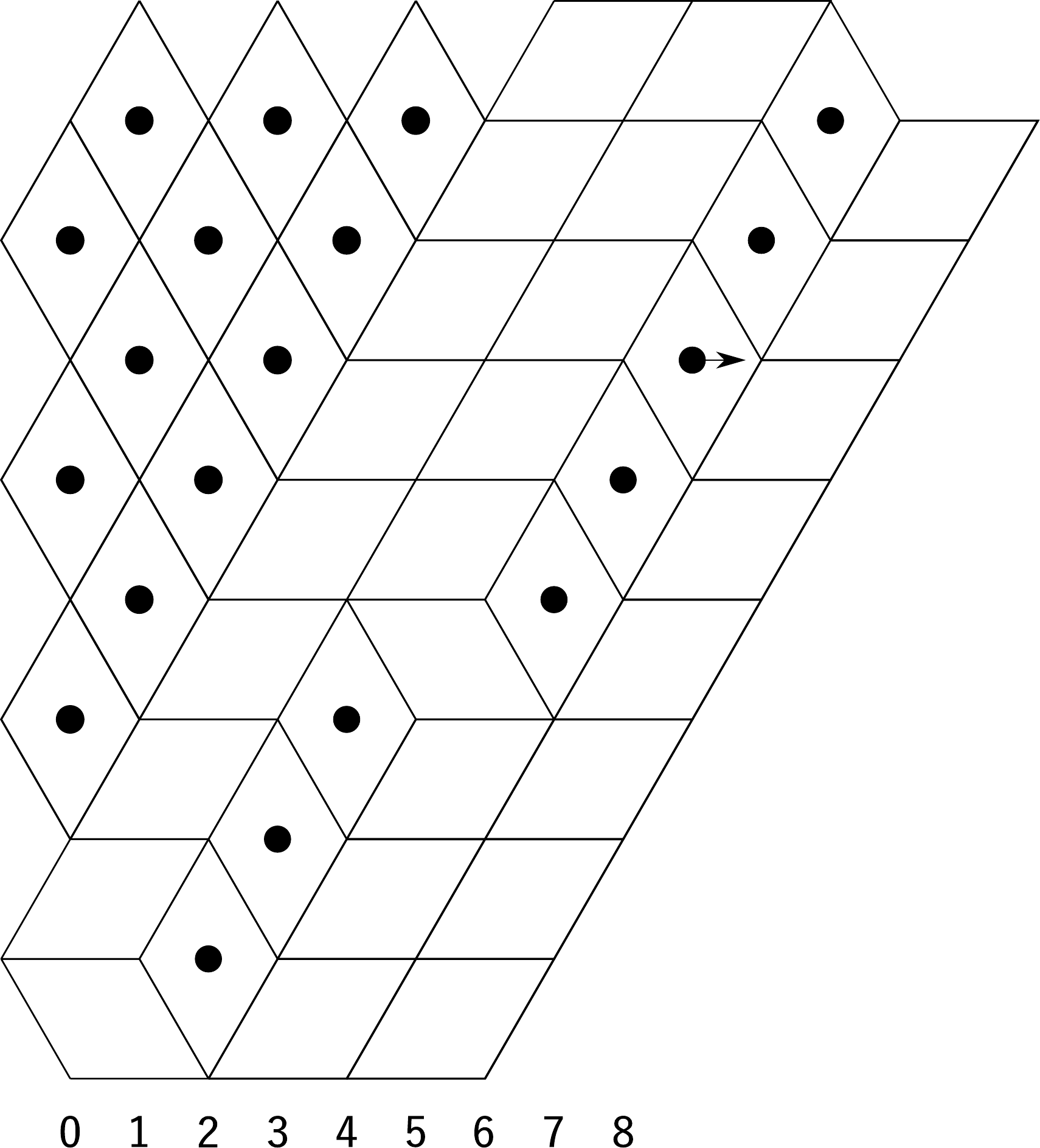}
\includegraphics[width=0.24\textwidth]{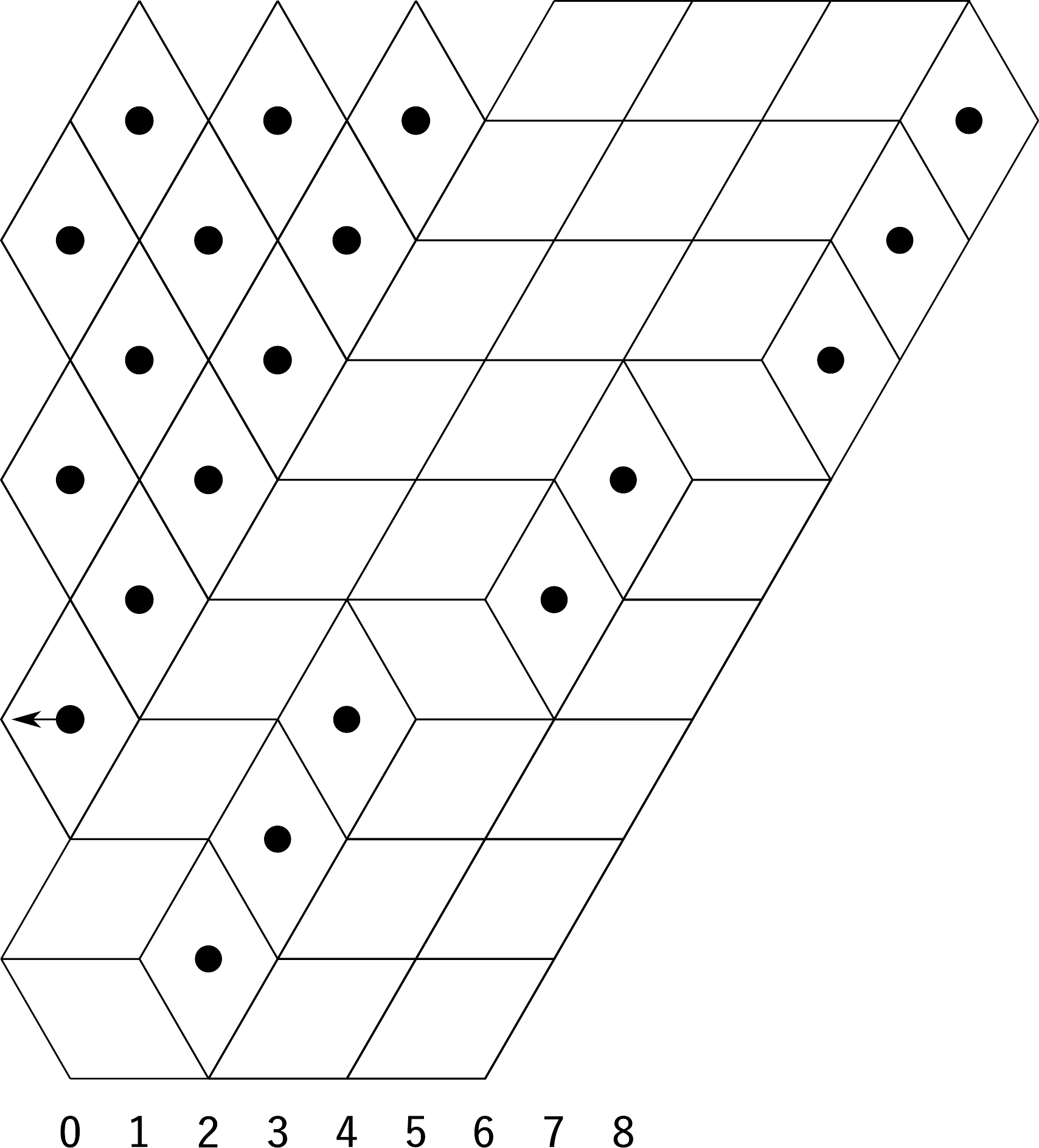} 
\caption{Three possible initial steps of either orthogonal or {symplectic} Plancherel growth. The arrows indicate the direction that the particle will attempt to jump next; see Figure \ref{modelcomparison} for the fourth step. The tiling determined by their positions suggests a stepped-surface interpretation.} \label{dynamicsteps}
\end{figure}

We now view either system of interacting particles as determining a random subset $\mc{X}(\gamma) := \{ \mc{X}^n_k(\gamma) \}_{\substack{ n \geq 1 \\ 1 \leq k \leq r_n}}$ of $\Z_{\geq 0} \times \Z_{> 0}$. We write $\xi^\gamma$ for the distribution on $2^{\Z_{\geq 0} \times \Z_{> 0}}$ of such point configurations {\color{black}(we refer to this distribution as a \emph{point process})}. Borodin--Kuan \cite{abjk1} studied the local asymptotics of such systems and their results are summarized in Figure \ref{averagebehavior}. Their arguments rely on the fundamental result that this point process is determinantal with an explicitly computable kernel; take $\alpha=\beta=-1/2$ in {\color{black} Theorem \ref{kernel}} for a precise statement. As indicated by this statement, certain orthogonal polynomials underlie the correlation structure, and the dynamics given by driving $\mbf{simple}$ random walks follow from the recurrence
\begin{equation} \label{chebyshevrecurrence}
xp_n(x) = \frac{1}{2} [p_{n+1}(x) + p_{n-1}(x)], \quad{\color{black} n>0}
\end{equation}
satisfied by the relevant first/fourth (orthogonal case) and second/third ({symplectic} case) Chebyshev polynomials (see e.g. \cite{jcm} for the definition of the third and fourth Chebyshev polynomials). Indeed, these are the only four Jacobi polynomials $P^{(\alpha, \beta)}_k$ that satisfy the recurrence \eqref{chebyshevrecurrence} (after an appropriate scaling -- {\color{black} see \eqref{recurrence} below}). The cases $\beta =\pm 1/2$ of the multi-time versions of the kernels  arise from the asymptotics at $-1$ of these polynomials.

Having reviewed these two archetypes, we can quickly introduce the \emph{Jacobi growth process}. Fix $\alpha, \beta>-1$. Consider the evolving point configuration $\mc{X}_{(\alpha,\beta)}(\gamma) := \{ (\mc{X}_{(\alpha,\beta)})^n_k(\gamma) \}_{\substack{ n \geq 1 \\ 1 \leq k \leq r_n}}$ of particles in $\Z_{\geq 0} \times \Z_{> 0}$ with the same point pattern preserved by push--block interactions as in the orthogonal and {symplectic} cases, except we now declare that if a particle is at position $\lambda \geq 0$ on level $n \geq 1$ and is allowed to jump to position $\mu$, then it will do so at the rate
\begin{equation} \label{generalrate}
r_n(\lambda, \mu):=
\begin{cases}
\begin{cases}
\frac{\lambda+\beta+1}{2\lambda+\alpha+\beta+2}  \frac{2(\lambda+1)}{2\lambda+\alpha+\beta+3}, & n \ \ \text{even} \\
\frac{\lambda+\alpha+\beta+1}{2\lambda+\alpha+\beta+1}  \frac{2(\lambda+\alpha+1)}{2\lambda+\alpha+\beta+2}, & n \ \ \text{odd}
\end{cases}, & \lambda \geq 0, \mu = \lambda+1 \\
\begin{cases}
\frac{\lambda+\alpha+\beta+1}{2\lambda+\alpha+\beta+1}  \frac{2(\lambda+\alpha+1)}{2\lambda+\alpha+\beta+2}, & n \ \ \text{even} \\
\frac{\lambda+\beta}{2\lambda+\alpha+\beta}  \frac{2\lambda}{2\lambda+\alpha+\beta+1}, & n \ \ \text{odd}
\end{cases}, & \lambda \geq 1, \mu = \lambda -1 
\end{cases},
\end{equation}
where $r_n(0,-1) \equiv 0$. {\color{black} The initial condition is once again the {\color{black} {densely packed} initial condition $\mc{X}^n_k(0) = r_n-k$}.} The reader can check that the particle system corresponds to the orthogonal case if $\alpha = \beta = -1/2$ and to the {symplectic} case if $-\alpha = \beta = 1/2$; in particular, the leftmost particles on odd levels have rate $\frac{2(\alpha+1)}{(\alpha+\beta+2)}$ right jumps off of the wall at $0$. { Also note that the jump rates converge to $1/2$ as $\lambda\rightarrow\infty$}. These rates have a more enlightening and compact form in terms of quantities arising from Jacobi polynomials (see equation \eqref{ratejacobiform} below). 

{\color{black}
It is worth mentioning that there are several other asymptotic regimes that could be pursued. 
In addition to investigating the edge asymptotic behavior, \cite{abjk1} also identifies the incomplete beta kernel (defined in \cite{aonr1}) as describing the asymptotics in the bulk of the system. Some formal calculations based on the results of Chapter 8 of \cite{gs1} suggest the same is true for the Jacobi growth process, but we do not pursue the rigorous details here. These formal calculations should also imply Gaussian free field fluctuations in the height function, using the results of \cite{jkgff}.

The paper \cite{BFPSW} identifies and \cite{cerenzia1} constructs the limit of orthogonal and symplectic Plancherel growth under diffusive scaling (i.e., the bottom left corner of Figure [3]). The fact that the drift of the rates (4) of the driving random walks vanishes suggests the Jacobi growth process shares the same diffusion limit. See \cite{acw} for recent work on such systems.
}

\section{Jacobi Polynomials}
{\color{black} In order to state the main results of this paper, let us first introduce some notation concerning the Jacobi polynomials.}

Define the Jacobi Polynomials $\{P^{(\alpha,\beta)}_k \}_{k \geq 0}$ by the initial condition $P^{(\alpha,\beta)}_0(x) \equiv 1$ and by the recurrence 
\begin{equation} \label{recurrence}
\ba
&x P^{(\alpha, \beta)}_k(x) = A_k^{(\alpha,\beta)} P^{(\alpha, \beta)}_k(x)+  B_k^{(\alpha,\beta)}  P^{(\alpha, \beta)}_{k-1}(x) + C_k^{(\alpha,\beta)}  P^{(\alpha, \beta)}_{k+1}(x), \ \ \ k \geq 0.
\ea
\end{equation}
where $P^{(\alpha,\beta)}_{-1} \equiv 0$ and where
$$
A_k^{(\alpha,\beta)} :=  \frac{(\beta - \alpha)(\alpha+\beta)}{(2k + \alpha + \beta)(2k + \alpha + \beta+2)} \ \ \ 
C_k^{(\alpha,\beta)} := \frac{2(k+1)(k+1+\alpha + \beta)}{(2k + \alpha + \beta+1)(2k + \alpha + \beta + 2)}, \ \ k \geq 0 
$$
$$
B_k^{(\alpha,\beta)} :=\frac{2 (k+\alpha)(k+\beta)}{(2k + \alpha + \beta)(2k + \alpha + \beta+1)}, \ \ k \geq 1.
$$
{\color{black} Observe that for any values of $\alpha$ and $\beta$, 
$$
\lim_{k\rightarrow\infty} A_k^{(\alpha,\beta)} = 0, \quad \quad \lim_{k \rightarrow \infty} B_k^{(\alpha,\beta)} = \lim_{k \rightarrow\infty} C_k^{(\alpha,\beta)}=\frac{1}{2}.
$$
}

The $P^{(\alpha,\beta)}_k$ are orthogonal with respect to the weight $w_{(\alpha,\beta)}(x):=(1-x)^\alpha (1+x)^\beta 1_{[-1,1]}(x)$, with normalization $P^{(\alpha,\beta)}_k(1) = \frac{\Gamma(k+\alpha+1)}{\Gamma(k+1) \Gamma(\alpha+1)}$, so that the leading coefficient is given by
\begin{equation} \label{leadingcoef}
a^{(\alpha,\beta)}_k : = \frac{\Gamma(2k+\alpha+\beta+1)}{2^k k! \Gamma(k+\alpha+\beta+1)} = \prod_{\ell=0}^{k-1} (C_\ell^{(\alpha,\beta)})^{-1}
\end{equation}
and the inner product by
\begin{equation} \label{innerproduct}
\langle P^{(\alpha,\beta)}_k, P^{(\alpha,\beta)}_\ell \rangle_{(\alpha,\beta)}:= \int_\R P^{(\alpha,\beta)}_k(x)P^{(\alpha,\beta)}_\ell(x) w_{(\alpha,\beta)}(x) dx = \frac{2^{\alpha+\beta+1} \Gamma(k+\alpha+1) \Gamma(k+\beta+1)}{(2k+\alpha+\beta+1)k!\Gamma(k+\alpha+\beta+1)} \delta_{k\ell}.
\end{equation}
We will always write $ \tilde{P}^{{\color{black} (\alpha, \beta)}}_k :=  P^{(\alpha,\beta)}_k/ \Vert P^{(\alpha,\beta)}_k\Vert$ so that $\langle  \tilde{P}^{{\color{black} (\alpha, \beta)}}_k,  \tilde{P}^{{\color{black} (\alpha, \beta)}}_\ell \rangle_{(\alpha, \beta)}= \delta_{k\ell}$. We will frequently make use of the orthogonal decomposition
\begin{equation} \label{orthogonaldecomposition}
T(x) = \sum_{k=0}^\infty \langle  \tilde{P}^{{\color{black} (\alpha, \beta)}}_k, T \rangle_{(\alpha, \beta)}  \tilde{P}^{{\color{black} (\alpha, \beta)}}_k(x)
\end{equation} for $T \in C^1[-1,1]$. See Szeg{\color{black} \H{o}} \cite{gs1} for more details on this discussion.

Fix $\alpha, \beta >-1$; the parameter $\beta$ will hereafter be suppressed from notation whenever possible. Write
$$
\alpha_n:=
\begin{cases}
\alpha+1, & \text{if} \ n \ \text{even} \\
\alpha, & \text{if} \ n \ \text{odd}
\end{cases}.
$$
{\color{black} The use of this $ \alpha_n$ has occurred before, see (7.4) of \cite{aogo1}. }
Besides the orthonormal scaling $\tilde{P}_k^{(\alpha)}$, it will be convenient to define\footnote{ {\color{black} Note that there is an abuse of notation here, in that the scalings $\bar{P}_k^{(\alpha_n )}$ and $\bar{\bar{P}}_k^{(\alpha_n)}$ depend on the parity of $n$. So for example, $\bar{P}_k^{(\alpha+1)}$ should be interpreted as $\bar{c}_k^2 P_k^{(\alpha+1)}$ and $\bar{P}_k^{(\alpha)}$ should be interpreted as $\bar{c}_k^1 P_k^{(\alpha)}$.  } }
$$\bar{P}^{(\alpha_n)}_k = \bar{c}_k^n P^{(\alpha_n)}_k, \ \ \ \ \bar{\bar{P}}^{(\alpha_n)}_k =\bar{\bar{c}}_k^n P^{(\alpha_n)}_k,$$ where
\begin{equation}  \label{scalings}
 \bar{c}_k^n : = (2k+\alpha_n+\beta + 1) \cdot
\begin{cases}
\frac{\Gamma(k+1) \Gamma(\alpha+1)}{2\cdot \Gamma(\alpha+k+2)}, & n \ even \\
\frac{\Gamma(k+\alpha+\beta + 1)}{\Gamma(k+\beta+1)}, & n \ odd
\end{cases} \ \ \  
\bar{\bar{c}}_k^n : = 
\begin{cases}
\frac{\Gamma(k+\alpha+\beta+2)}{\Gamma(k+\beta+1)}, & n \ even \\
\frac{\Gamma(k+1) \Gamma(\alpha+1)}{\Gamma(\alpha+k+1)}, & n \ odd
\end{cases}
\end{equation}
Define
$$
\phi_{n}(k):= \bar{c}^{n}_k/\bar{\bar{c}}^{n}_k, \ \ \  \ \ 
\phi_{n}(k,m) :=
\phi_n(k) \cdot
\begin{cases}
1_{(k<m)},&\text{if} \ n \ \text{even} \\ 
1_{(k \leq m)},& \text{if} \ n \ \text{odd}
\end{cases}.
$$ 
We also adopt the convention $\phi_n(-1, \cdot) \equiv 1$. Note the facts
\begin{equation} \label{barsidentity}
\phi_{n}(s) \bar{\bar{P}}_s^{(\alpha_{n})} = \bar{P}_s^{(\alpha_{n})}, \ \ \ \bar{\bar{P}}^{(\alpha)}_k(1) = 1,
\end{equation}
and that the scalings are related by
\begin{equation} \label{scaling}
\frac{\bar{P}_k^{(\alpha_n)}\bar{\bar{P}}_k^{(\alpha_n)}}{2^{\alpha+\beta+1} \Gamma(\alpha+1)} = \tilde{P}_k^{(\alpha_n)} \tilde{P}_k^{(\alpha_n)}.
\end{equation}
Finally, with this notation established, we can rewrite the rate \eqref{generalrate} compactly as
\begin{equation} \label{ratejacobiform}
r_n(\lambda, \mu)=
\begin{cases}
\frac{\bar{c}_{\lambda+1}^n}{\bar{c}_{\lambda}^n}  B^{(\alpha_n,\beta)}_{\lambda+1} = \frac{\bar{\bar{c}}_{\lambda}^n}{\bar{\bar{c}}_{\lambda+1}^n} C^{(\alpha_n,\beta)}_{\lambda}, & \lambda \geq 0, \mu = \lambda+1 \\
\frac{\bar{c}_{\lambda-1}^n}{\bar{c}_{\lambda}^n} C^{(\alpha_n,\beta)}_{\lambda-1} = \frac{\bar{\bar{c}}^n_{\lambda}}{\bar{\bar{c}}^n_{\lambda-1}}  B^{(\alpha_n,\beta)}_{\lambda}, & \lambda \geq 1, \mu = \lambda -1 \\
\end{cases}.
\end{equation}

\section{Statement of main theorems}

Below, we write $\xi^\gamma = \xi^\gamma _{(\alpha,\beta)}$ for the distribution on $2^{\Z_{\geq 0} \times \Z_{>0}}$ of the random point configuration $\mc{X}_{(\alpha,\beta)}(\gamma) \in 2^{\Z_{\geq 0} \times \Z_{>0}}$ determined by the Jacobi growth dynamics.

\bt \label{kernel}
Fix $\alpha, \beta >-1$ and $\gamma \geq 0$. Then the correlation functions $\{ \rho^\gamma_k\}_{k \geq 1}$ of the Jacobi growth process are determinantal: for $z_1, \ldots, z_k \in \Z_{\geq 0} \times \Z_{> 0}$,
$$
\rho^\gamma_k(z_1,\ldots,z_k) := \xi^\gamma(\{ E \in 2^{\Z_{\geq 0} \times \Z_{> 0}} | E \supset \{ z_1, \ldots, z_k \} \})= \det[ K^\gamma_{(\alpha,\beta)}(z_i, z_j) ]_{i,j=1}^k,
$$
where the (nonsymmetric) kernel $K^\gamma_{(\alpha,\beta)}$ is given by
\begin{equation} \label{explicitkernel}
\ba
K^\gamma_{(\alpha,\beta)}((s,n),(t,m)) :=    & \frac{1}{2 \pi i} \int_{-1}^1 \oint \frac{e^{x\gamma}}{e^{u\gamma}}\frac{\bar{P}_{s}^{(\alpha_n)}(x)}{2^{\alpha+\beta+1}\Gamma(\alpha+1)}  \bar{\bar{P}}_{t}^{(\alpha_m)}(u) \frac{(1-x)^{r_n + \alpha_n}(1+x)^{\beta}}{(1-u)^{r_m}(x-u)} dudx\\
& + 1_{(n \geq m)} \int_{-1}^1 \frac{\bar{P}_{s}^{(\alpha_n)}(x)}{2^{\alpha+\beta+1}\Gamma(\alpha+1)} \bar{\bar{P}}_{t}^{(\alpha_m)}(x) (1-x)^{r_n-r_m+\alpha_n} (1+x)^{\beta} dx. 
\ea
\end{equation}
for $(s,n),(t,m) \in \Z_{\geq 0} \times \Z_{> 0}$, where the complex integral is a positively oriented (i.e., counterclockwise) simple loop around $[-1,1]$.
\et

At the edges of the particle system, we find an $(\alpha,\beta)$-dependent family of Jacobi kernels and a $\beta$-dependent family of {\color{black} hard--edge} Pearcey kernels. We first focus on the large time limit at a finite distance from the wall. 

\bt{(Edge Limit: Discrete Jacobi Kernel)} \label{jacobilimit}
Fix $\alpha, \beta >-1$. Assume $\gamma \geq 0$, $(s_i, n_i) \in \Z_{\geq 0} \times \Z_{> 0}$, $1 \leq i \leq k$, depend on $N$ in such a way that $\gamma \sim N  \tau >0$, $r_{n_i} \sim N  \eta>0 $ but the $s_i$ are fixed and finite. Assume the parity of each $n_i$ is constant and that the differences $r_{n_i} - r_{n_j}$ are of constant order. Then as $N \to \infty$,
$$
\rho^\gamma_k  ((s_1, n_1), \ldots ,(s_k, n_k)) \to
\begin{cases}
\det [ \mc{J}_{(\alpha,\beta)}((s_i, n_i), (s_j, n_j)|1-\eta/\tau) ]_{i,j=1}^k, & \text{if} \ \ 1-\eta/\tau > -1 \\
1, & \text{if} \ \ 1-\eta/\tau \leq -1
\end{cases}
$$
where the discrete Jacobi kernel $\mc{J}_{(\alpha,\beta)}$ is given by
\begin{equation} \label{explicitjacobi}
\ba
\mc{J}_{(\alpha,\beta)}((s,n), (t,m) | \epsilon) : = & \int_{-1}^1 [1_{(n\geq m)}+ 1_{[-1,\epsilon]}(x)]  \frac{\bar{P}_{s}^{(\alpha_n)}(x)}{2^{\alpha+\beta+1}\Gamma(\alpha+1)}   \bar{\bar{P}}_{t}^{(\alpha_m)}(x) (1-x)^{r_n-r_m+\alpha_n} (1+x)^{\beta} dx \\
\ea.
\end{equation}
\et
Define the \emph{particle--hole involution} as the map $\Delta$ on $2^{\Z_{\geq 0} \times \Z_{>0}}$ given by $\Delta(\mc{X}) := (\Z_{ \geq 0} \times \Z_{> 0}) \setminus \mc{X}$. The push--forward $\xi^\gamma_\Delta : = \xi^\gamma \circ \Delta^{-1}$ of the point process $\xi^\gamma$ furnishes another point process which also possesses a determinantal correlation structure; see Proposition 3.1 of \cite{abjk2} and equation \eqref{compkernel}. Our main result shows that the kernel arising in our setting at the critical point {\color{black} is the hard--edge Pearcey kernel from \cite{sdbv1}}. Note in particular the result is independent of $\alpha>-1$. Below, recall that $J_{\beta}$ is the Bessel function of the first kind, {\color{black} and $I_{\beta}$ is the modified Bessel function of the first kind, defined by} 
$$
J_\beta(x) = \sum_{k=0}^\infty \frac{(-1)^k (x/2)^{\beta+2k}}{k! \Gamma(\beta+k+1)}, \quad {\color{black} I_\beta(x) = \sum_{k=0}^\infty \frac{ (x/2)^{\beta+2k}}{k! \Gamma(\beta+k+1)}}.
$$ 
{\color{black}
The hard--edge Pearcey kernel $L^{\beta}$ of \cite{sdbv1} is defined by 
\begin{multline*}
\left( \frac{x}{y} \right)^{\beta/2} L^{\beta}(s,x;t,y) = - 1_{\{t>s\}} \frac{1}{t-s} \exp\left( - \frac{x+y}{t-s} \right) I_{\beta}\left( \frac{ 2\sqrt{xy}}{t-s}\right) \\
 +  \frac{2}{\pi i} \int_C dv \int_0^{\infty} du \left( \frac{u}{v} \right)^{\beta} \frac{uv}{v^2-u^2} \frac{e^{v^4/2+sv^2}}{e^{u^4/2+tu^2}} J_{\beta}(2\sqrt{ y} u)J_{\beta}(2 \sqrt{x} v),
\end{multline*}
where $C$ consists of two rays, one from $e^{i \pi /4} \infty$ to $0$ and one from $0$ to $e^{-i \pi / 4}\infty$. Note that the $(x/y)^{\beta/2}$ term is a conjugating factor which disappears in the determinant.}

\bt{(Edge Limit: {\color{black} Hard--edge} Pearcey Kernel)} \label{pearceylimit}
Fix $\alpha, \beta > -1$. Let $\rho^{\gamma, k}_\Delta$ denote the $k$th correlation function of the point process $\xi^\gamma_\Delta$: for distinct $z_1, \ldots , z_k \in \Z_{\geq 0} \times \Z_{> 0}$, 
$$
\rho^{\gamma, k}_\Delta(z_1, \ldots, z_k) : = \xi^\gamma_\Delta(\{ E \in 2^{\Z_{\geq 0} \times \Z_{> 0}} | E \supset \{ z_1, \ldots, z_k \} \}).
$$
Assume $\gamma \geq 0$, $(s_i, n_i) \in \Z_{\geq 0} \times \Z_{> 0}$, $1 \leq i \leq k$, depend on $N$ in such a way that $\gamma \sim N/2$, $s_i \sim N^{1/4} \nu_i>0$, and $r_{n_i} -N\sim \sqrt{N} \sigma_i$. Assume the parity of each $n_i$ is constant. Then as $N \to \infty$
$$
\ba
 N^{k/4}  \cdot \rho^{\gamma, k}_\Delta & ((s_1, n_1), \ldots ,(s_k, n_k))
\to
\det [ {\color{black} L^{\beta}(\sigma_i,\nu_i^2;\sigma_j,\nu_j^2) \cdot 2\sqrt{\nu_i\nu_j} } ]_{i,j=1}^k.
\ea
$$
\et
 
{\color{black} 
 
\begin{remark} The expression $2\sqrt{\nu_i\nu_j}$ arises due to the change of variables $x = \nu_i^2, y= \nu_j^2$ in the spatial dimensions. This change of variables is related to the fact that $L^{\beta}$ occurred in the context of \textit{squared} Bessel paths. Indeed, it can be written as
$$
2\sqrt{\nu_i\nu_j} = \sqrt{ \frac{dx}{d \nu_i} \frac{dy}{d \nu_j}  }.
$$ 
See Remark \ref{CoV} below for a more in--depth explanation.
\end{remark}
}

{\color{black} Let us review the outline of the remainder of the paper. Section 3 will prove some necessary identities for the Jacobi polynomials. Section 4 constructs probability measures on the interlacing particle configurations, and section 5 constructs a Markov process preserving these measures. Section 6 computes an explicit formula for the correlation kernel, and section 7 takes the asymptotics of this kernel.

\textbf{Acknowledgements} The authors would like to thank Alexei Borodin for valuable discussion. M.C. was supported by NSF grant DMS--0806591 and J.K. was supported by the Minerva Foundation and NSF grant DMS--1502665.
}

\section{Some ancillary results} \label{sectionnotations}

We now collect some results that are fundamental for computing the correlation kernel in Section \ref{sectiondeterminantal} and for uncovering the intertwining relationship in Section \ref{sectionMarkov}. The first generalizes Lemma 2.4 and 2.5 of \cite{abjk1}. 

\bp  \label{identities}
The Jacobi polynomials satisfy the following identities: for any $T \in C^1[-1,1]$, 
\begin{enumerate}

\item

\subitem $\sum_{r=0}^{s}  \bar{P}^{(\alpha)}_r =  \bar{\bar{P}}^{(\alpha+1)}_s$, for all $s \geq 0$. 

\subitem $\sum_{r =s+1}^\infty \left \langle \bar{P}^{(\alpha)}_r, T \right \rangle_{\alpha}  = \langle \bar{\bar{P}}^{(\alpha+1)}_s, T(1)-T \rangle_{\alpha} $

\item 

\subitem $ \sum_{r=0}^{s-1} \bar{P}^{(\alpha+1)}_r = \frac{ \bar{\bar{P}}^{(\alpha)}_s- 1}{x-1}$, for all $s \geq 1$.

\subitem  $ \sum_{r=s}^\infty  \langle \bar{P}^{(\alpha+1)}_r , T \rangle_{\alpha+1} =  \langle \bar{\bar{P}}^{(\alpha)}_s , T \rangle_{\alpha}$

\end{enumerate}
\ep

\bpf

For the first identity, note that if $p(x)$ is any polynomial of degree $<s$, then 
$$
\left \langle \sum_{r=0}^{s}  \bar{P}^{(\alpha)}_r, p \right \rangle_{\alpha+1} = \sum_{r=0}^{s}  \left \langle  \bar{P}^{(\alpha)}_r, p \cdot (1-x) \right \rangle_{\alpha} \bar{\bar{P}}^{(\alpha)}_r(1) = 2^{\alpha+\beta+1} \Gamma(\alpha+1) \cdot p(x) (1-x) |_{x=1} = 0,
$$
where we have used the facts $\bar{\bar{P}}^{(\alpha)}_r(1) = 1$ and \eqref{scaling} along with orthogonal decomposition \eqref{orthogonaldecomposition}. Hence, by uniqueness of orthogonal polynomials, $\sum_{r=0}^{s} \bar{P}^{(\alpha)}_r$ must be proportional to $P^{(\alpha+1)}_s$ and comparing leading coefficients tells us that the proportionality constant is (see equations \eqref{leadingcoef} and \eqref{scalings})
$$
(2s+\alpha+\beta+1)\frac{\Gamma(s+\alpha+\beta + 1)}{\Gamma(s+\beta+1)} \cdot \frac{a^{(\alpha, \beta)}_s}{a^{(\alpha+1,\beta)}_s} = \frac{\Gamma(s+\alpha+\beta+2)}{\Gamma(s+\beta+1)},
$$
as required. On the other hand, combining this identity with the Christoffel--Darboux formula gives
\begin{equation} \label{CDformula}
\ba
\bar{\bar{P}}^{(\alpha+1)}_s(x) &= \sum_{r=0}^{s}  \bar{P}^{(\alpha)}_r(x) \bar{\bar{P}}^{(\alpha)}_r(1) = \frac{2^{\alpha+\beta+1} \Gamma(\alpha+1)  a^{(\alpha, \beta)}_s}{ \Vert P^{(\alpha,\beta)}_s \Vert^2 a^{(\alpha, \beta)}_{s+1}} \frac{P^{(\alpha,\beta)}_{s+1}(x) P^{(\alpha,\beta)}_s(1) - P^{(\alpha,\beta)}_s(x) P^{(\alpha,\beta)}_{s+1}(1)}{x-1} \\
& = \frac{2 \cdot \Gamma(\alpha+1) \Gamma(s+2) \Gamma(s+\alpha+\beta+2)}{(2s+\alpha+\beta+2)\Gamma(s+\alpha+1)\Gamma(s+\beta+1)} \cdot \frac{P^{(\alpha,\beta)}_{s+1}(x) P^{(\alpha,\beta)}_s(1) - P^{(\alpha,\beta)}_s(x) P^{(\alpha,\beta)}_{s+1}(1)}{x-1}
\ea
\end{equation}
Multiply both sides of equation \eqref{CDformula} by
$$
\phi_n(s) = \frac{\bar{c}^n_s}{\bar{\bar{c}}^n_s} = \frac{(2s+\alpha+\beta+2) \Gamma(s+1) \Gamma(\alpha+1) \Gamma(s+\beta+1)}{2 \cdot \Gamma(s+\alpha+2) \Gamma(s+\alpha+\beta+2)}, \  \ n \ \ \text{even}
$$
and recall 
$$
 P^{(\alpha,\beta)}_s(1) = \frac{\Gamma(s+\alpha+1)}{\Gamma(s+1) \Gamma(\alpha+1)}.
$$
Then the expression \eqref{CDformula} becomes
\begin{equation} \label{okounkov7.7}
\bar{P}^{(\alpha+1)}_r =  \frac{\bar{\bar{P}}^{(\alpha,\beta)}_{r+1}(x) - \bar{\bar{P}}^{(\alpha,\beta)}_{r}(x)}{x-1},
\end{equation} 
which telescopes when summing over $r$ to give the first identity of the second point.

For the latter identities of each point, first note that since $P^{(\alpha)}_0 \equiv 1$, 
$$
\ba
\langle \bar{\bar{P}}^{(\alpha+1)}_s , 1 \rangle_{\alpha} & = \left \langle \sum_{r=0}^{s}  \bar{P}^{(\alpha)}_r , P^{(\alpha)}_0 \right \rangle_{\alpha} = \frac{(\beta +\alpha+1) \Gamma(\alpha+\beta+1)}{\Gamma(\beta+1)} \langle P^{(\alpha)}_0, P^{(\alpha)}_0 \rangle_{\alpha} =2^{\alpha+\beta+1}\Gamma(\alpha+1)
\ea
$$
This fact along with orthogonal decomposition \eqref{orthogonaldecomposition} gives
$$
\ba
\langle \bar{\bar{P}}^{(\alpha+1)}_s, T(1) \rangle_{\alpha} & = 2^{\alpha+\beta+1}\Gamma(\alpha+1) T(1) = 2^{\alpha+\beta+1}\Gamma(\alpha+1) \sum_{r =0}^\infty \langle \tilde{P}^{(\alpha)}_r, T \rangle_{\alpha} \tilde{P}^{(\alpha)}_r(1) = \sum_{r =0}^\infty \left \langle \bar{P}^{(\alpha)}_r, T \right \rangle_{\alpha} 
 \ea
$$
Then subtract $ \langle \bar{\bar{P}}^{(\alpha+1)}_s, T \rangle_{\alpha}  = \left \langle \sum_{r=0}^s \bar{P}^{(\alpha)}_s, T \right \rangle_{\alpha} $ from both sides to arrive at the second identity of the first item. For the last identity, note
$$
\left \langle \sum_{r=0}^{s-1} \bar{P}^{(\alpha+1)}_r, T \right \rangle_{\alpha+1}= \left \langle \frac{1-  \bar{\bar{P}}^{(\alpha)}_s}{1-x}, T \right \rangle_{\alpha+1} \overset{s \to \infty}{\rightarrow} \langle 1 , T \rangle_{\alpha}
$$
where the term disappears necessarily from convergence of orthogonal decompositions \eqref{orthogonaldecomposition}. Subtracting the finite sum $\sum_{r=0}^{s-1} \langle \bar{P}^{(\alpha+1)}_r, T \rangle_{\alpha+1} = \langle 1- \bar{\bar{P}}^{(\alpha)}_s, T \rangle_{\alpha} $ from the right side completes the proof.

\epf
For smooth $E \in C^\infty[-1,1]$, denote the $m$th Taylor Remainder of $E$ about $1$ by
$$
R^E_m(x): = 
\begin{cases}
E(x), & \ \ m \leq 0 \\
E(x) - \sum_{k=0}^{m-1} \frac{E^{(k)}(1)}{k!}(x-1)^{k}, & \ \ m \geq 1 
\end{cases}.
$$
and let 
$$
\Psi^{n}_{r_n-l}(s) = \Psi^{n, E}_{r_n-l}(s) : =  \left  \langle \frac{\bar{P}^{(\alpha_n)}_s}{2^{\alpha+\beta+1}\Gamma(\alpha+1)},  (x-1)^{r_n-l} R^{E}_{l-r_n} \right \rangle_{\alpha+1}, 
$$
The next result generalizes Theorem 4.4 in \cite{abjk1}. 
\bl \label{comprule}
The functions $\Psi^{n}_{r_n-l}(s)$, $n \geq 1$, $l \in \Z$, satisfy the composition rule
$$
(\phi_{n-1}*\Psi^{n}_{r_n-l})(s):= \sum_{r\geq 0} \phi_{n-1}(s,r)\Psi^{n}_{r_n-l}(r) = \Psi^{n-1}_{r_{n-1}-l}(s).
$$
\el
\bpf
Fix $n$ odd so $r_{n-1} = r_n-1$ and $\alpha_{n} = \alpha$. Then using part $(3)$ of Proposition \ref{identities} with $T_m(x) := (x-1)^{-m} R^E_{m}(x)$ and the fact $T_m(1) := \lim_{x \to 1} T_m(x) = E^{(m)}(1)/m!$, we get using \eqref{barsidentity}
$$
\ba
\sum_{r\geq 0} \phi_{n-1}(s,r)\Psi^{n}_{r_n-l}(r) & = \frac{\phi_{n-1}(s)}{2^{\alpha+\beta+1}\Gamma(\alpha+1)}  \sum_{r=s+1}^\infty \langle \bar{P}^{(\alpha)}_r, (x-1)^{r_n-l} R^E_{l-r_n} \rangle_{\alpha} \\
&= \frac{1}{2^{\alpha+\beta+1}\Gamma(\alpha+1)} \left \langle \phi_{n-1}(s) \bar{\bar{P}}^{(\alpha+1)}_s, \frac{T_{r_n-l}-T_{r_n-l}(1)}{x-1} \right \rangle_{\alpha+1} \\
&  = \frac{1}{2^{\alpha+\beta+1}\Gamma(\alpha+1)}  \left \langle \bar{P}^{(\alpha+1)}_s, \frac{R^E_{l-r_{n-1}}(x)}{(x-1)^{l-r_{n-1}}} \right \rangle_{\alpha+1} = \Psi^{n-1}_{r_{n-1}-l}(s),
\ea
$$
The case $n$ even instead involves $r_n = r_{n-1}$, $\alpha_{n} = \alpha+1$, and the other part of Proposition \ref{identities}:
$$
\ba
\sum_{r\geq 0} \phi_{n-1}(s,r)\Psi^{n}_{r_n-l}(r)& = \frac{1}{2^{\alpha+\beta+1}\Gamma(\alpha+1)}  \sum_{r=s}^\infty \langle  \bar{P}^{(\alpha+1)}_r, (x-1)^{r_n-l} R^E_{l-r_n} \rangle_{\alpha+1} \\
&=\frac{1}{2^{\alpha+\beta+1}\Gamma(\alpha+1)}  \langle \phi_{n-1}(s) \bar{\bar{P}}^{(\alpha)}_s, (x-1)^{r_{n-1}-l} R^E_{l-r_{n-1}} \rangle_{\alpha} = \Psi^{n-1}_{r_{n-1}-l}(s).
\ea
$$ 
\epf

\section{Consistent series of probability measures on partitions}

For $n \geq 1$, we write $r_n:= \lfloor (n+1)/2 \rfloor$. Let $\mathbb{J}_n$ denote the set of \emph{partitions} $\lambda=( \lambda_1 \geq \ldots \geq \lambda_{r_n} \geq 0)$ of nonnegative integers of length at most $r_n$ (the \emph{length} of a partition $\lambda$ is the number of nonzero terms). For $\lambda \in \mathbb{J}_n, \mu \in \mathbb{J}_{n+1}$, write $ \lambda \prec \mu$ if $\lambda$ \emph{interlaces} $\mu$ in the sense
$$
\mu_1 \geq \lambda_1 \geq \mu_2 \geq \lambda_2 \geq \ldots  \geq \mu_{r_n} \geq \lambda_{r_n} \geq \mu_{r_n+1}.
$$ 
For $\lambda \in \mb{J}_n$, the transformation $\widetilde{\lambda}_i : = \lambda_i +r_n - i$ arises naturally and frequently in our formulas. For the construction of Section \ref{sectioncentral}, define the collection $\mathbb{J}_{n,seq}$ of all finite sequences $u = (u^1, \cdots, u^n)$ of partitions $u^i \in \mb{J}_i$, $1 \leq i \leq n$, of length $n$. Let $\mb{J}_{n,paths} \subset \mb{J}_{n,seq}$ denote the subset of \emph{Gelfand-Tsetlin patterns}, i.e., finite interlaced sequences $u = (u^1 \prec \cdots \prec u^n)$. Let $\mb{J}_{\infty, seq} $, $\mb{J}_{\infty, paths} $ denote infinite versions of each of these sets. We adopt the following notational convention:

\begin{convention} 
Quantities with indices that overflow (e.g. $\lambda_{r_n +1}$ for $\lambda \in \mb{J}_n$) are set to $0$ and those with indices that underflow (e.g. $\lambda_{0}$ for $\lambda \in \mb{J}_n$) are set to $\infty$. 
\end{convention}

We now collect two elementary results that will allow us to express our probability measures on partition sequences in determinantal form and to work with the resulting expressions. The proof of the ifrst, which is straightforward and now standard, follows along the same lines as Proposition 3.5 of \cite{abjk2}, so we only state the version we need. 
\bl \label{indication}
Indication of interlacement between two partitions $\lambda \in \mb{J}_{n}, \mu \in \mb{J}_{n+1}$ takes the determinantal form
$$
1_{(\lambda \prec \mu)} =  
\begin{cases}
\det[ 1_{(\widetilde{\lambda}_i < \widetilde{\mu}_j)}]_{i,j=1}^{r_{n+1}}, & \text{if} \ \ n \ \ \text{even} \\
\det[ 1_{(\widetilde{\lambda}_i \leq \widetilde{\mu}_j)}]_{i,j=1}^{r_{n+1}}, & \text{if} \ \ n \ \ \text{odd}
\end{cases},
$$
where $\widetilde{\lambda}_i : = \lambda_i +r_n - i$.
\el
For reference, we recite a variant of the Cauchy--Binet formula.
\bl{(Lemma 2.1, \cite{abjk1})} \label{cauchybinet}
For each $k \geq 0$, let $f_k, g_k$ be functions on $\C$ such that $\sum_{k=0}^\infty f_k(w_i) g_k(z_j)$ converges absolutely for some $w_i, z_j \in \C$, $1 \leq i,j \leq n$. Then
$$
\sum_{k_1 > \ldots > k_n \geq 0} \det [ f_{k_j}(w_i)]_{i,j=1}^n \det [ g_{k_i}(z_j)]_{i,j=1}^n = \det \left [ \sum_{k=0}^\infty f_k(w_i) g_k(z_j)\right]_{i,j=1}^n.
$$
\el
We now introduce the fundamental probability measures on partitions and their cotransitions. For $\lambda \in \mb{J}_n$, let 
$$
{\color{black}
d_{\lambda}^{(n)}(x_1,\ldots,x_{r_n}) = d^{(\alpha;n,\beta)}_{\lambda}(x_1,\ldots,x_{r_n}) : = 
}\frac{\det \left [\bar{\bar{P}}_{\widetilde{\lambda}_i}^{(\alpha_n)}(x_j) \right ]_{i,j=1}^{r_n}}{\det [ x_j^{r_n-i} ]_{i,j=1}^{r_n}}.
$$
{ Note that this is well defined at $x_i=x_j$ because the denominator is a Vandermonde determinant, which divides the numerator. Up to a constant, $d_{\lambda}^{(n)}$ is a $BC_{r_n}$ orthogonal polynomial at $\theta=1$, see Proposition 7.1 of  \cite{aogo1}.}
\bp
The quantities
\begin{equation} \label{cotransition}
T^n_{n-1}(\lambda, \mu) :=  \prod_{k=1}^{r_{n-1}} \phi_{n-1}(\tilde{\mu}_k) 1_{(\mu \prec \lambda)} \frac{{\color{black} d_{\mu}^{(n-1)}(1,\ldots,1)}}{{\color{black} d_{\lambda}^{(n)}(1,\ldots,1) }} = \det[ \phi_{n-1}( \tilde{\mu}_i, \tilde{\lambda}_j) ]_{i,j =1}^{r_{n}}\frac{{\color{black} d_{\mu}^{(n-1)}(1,\ldots,1)}}{{\color{black} d_{\lambda}^{(n)}(1,\ldots,1)}},
\end{equation}
are \emph{cotransition probabilities} from $\lambda \in \mb{J}_{n}$ down to $ \mu \in \mb{J}_{n-1}$, $n \geq 2$. In particular, this furnishes stochastic operators $T^n_{n-1}$. 
\ep
\bpf
The second equality of \eqref{cotransition} follows from Lemma \ref{indication} and the convention $\phi_{n-1}(-1, \cdot) \equiv 1$ for $n$ odd. Stochasticity follows from a branching rule for ${\color{black} d_{\lambda}^{(n)}(1,\ldots,1)}$. If $n$ is odd so that $\alpha_n = \alpha$ and $r_{n-1} = r_n-1$, then
\begin{equation} \label{branchingrule}
\ba
\frac{\det \left [\bar{\bar{P}}_{\widetilde{\lambda}_i}^{(\alpha_n)}(x_j) \right ]_{i,j=1}^{r_n}}{\det [ x_j^{r_n-i} ]_{i,j=1}^{r_n}} \Bigg \vert_{x_{r_n} = 1} &= \frac{\det \left [\frac{\bar{\bar{P}}_{\widetilde{\lambda}_i}^{(\alpha)}(x_j) - \bar{\bar{P}}_{\widetilde{\lambda}_{i+1}}^{(\alpha)}(x_j)}{x_j-1} \right ]_{i,j=1}^{r_{n-1}}}{\det [ x_j^{r_{n-1}-i} ]_{i,j=1}^{r_{n-1}}} = \sum_{\mb{J}_{n-1} \ni \mu \prec \lambda} \frac{\det \left [\bar{P}_{\widetilde{\mu}_i}^{(\alpha+1)}(x_j) \right ]_{i,j=1}^{r_{n-1}}}{\det [ x_j^{r_{n-1}-i} ]_{i,j=1}^{r_{n-1}}} \\
&= \sum_{\mb{J}_{n-1} \ni \mu \prec \lambda}  \prod_{k=1}^{r_{n-1}} \phi_{n-1}(\tilde{\mu}_k) \frac{\det \left [\bar{\bar{P}}_{\widetilde{\mu}_i}^{(\alpha_{n-1})}(x_j) \right ]_{i,j=1}^{r_{n-1}}}{\det [ x_j^{r_{n-1}-i} ]_{i,j=1}^{r_{n-1}}},
\ea
\end{equation}
where we have used the fact $\bar{\bar{P}}_{\widetilde{\lambda}_i}^{(\alpha_n)}(1)=1$ and indicated properties of the determinants in the first equality, the identity \eqref{okounkov7.7} in the second, and the fact $\bar{P}_s^{(\alpha_{n})} = \phi_{n}(s) \bar{\bar{P}}_s^{(\alpha_{n})} $ in the last. If instead $n$ is even, then $r_n = r_{n-1}$ and the same identity holds without the evaluation ``$x_{r_n} = 1$" by using the first item in Proposition \ref{identities}.
\epf

\begin{remark} \label{combinatorialexpression}
Iterating the branching established in this proof shows the quantities ${\color{black} d_{\lambda}^{(n)}(1,\ldots,1)}$ have a combinatorial expression as a sum of weights over all paths of length $r_n$ ending at $\lambda$, which was already established as (7.5) and (7.6) in the proof of Proposition 7.5 of Okounkov--Olshanski \cite{aogo1}. The branching weights are built from the coefficients $B(m,\ell)$ of equation $(7.14)$ in \cite{aogo1}, which are given in our notation by
$$
B(m,\ell) = \phi_n(m) \phi_{n-1}(\ell) = \frac{\bar{c}^n_m}{\bar{\bar{c}}^n_m}\frac{\bar{c}^{n-1}_\ell}{\bar{\bar{c}}^{n-1}_\ell}, \ \ \ n \ \text{even}.
$$
\end{remark}

{  Let $\psi \in C^1[-1,1]$ }. The cotransition operators $T^n_{n-1}$ serve to link the series of measures on $\mb{J}_n$ defined by
\begin{equation} \label{plancherelmeasures}
\ba
P_n^{  { \psi} }(\lambda) :=
\det \left [ \frac{\langle \bar{P}^{(\alpha_n)}_{\widetilde{\lambda}_j} , (x-1)^{r_n-i}{ \psi} \rangle_{\alpha+1}}{2^{\alpha+\beta+1}\Gamma(\alpha+1)} \right]_{i,j=1}^{r_n} {\color{black} d_{\lambda}^{(n)}(1,\ldots,1) } = \det \left [ \Psi^{n,{ \psi} }_{r_n-i}(\widetilde{\lambda}_j) \right]_{i,j=1}^{r_n} {\color{black} d_{\lambda}^{(n)}(1,\ldots,1) } .
\ea
\end{equation}

\bp \label{partitionmeasure}
The series of measures $\{P^{ \psi}_n\}_{n=1}^\infty$ is consistent in the sense
$$
P^{{ \psi}}_{n-1}(\mu) = (P^{{ \psi}}_nT^n_{n-1})(\mu):= \sum_{\lambda \in \mb{J}_n} P^{{ \psi}}_n(\lambda) T^n_{n-1}(\lambda, \mu), \ \ \ n \geq 2,
$$
and each { satisfy
$$
\sum_{\mu \in \mb{J}_n} P^{\psi}_n(\mu) \frac{ d_{\mu}^{{\color{black}(n)}}(x_1,\ldots,x_{{\color{black} r_n}})}{d_{\mu}^{{\color{black} (n)}}(1,\ldots,1)} = \psi(x_1)\cdots \psi(x_{\color{black} r_n}),
$$
}
so in particular, have mass ${ \psi(1)^n}$. Moreover,  { if $\psi  \equiv 1$} then $P^{{ \psi}}_n(\lambda) = \delta_0(\lambda)$.
\ep

\bpf
For either parity of $n$, we deduce {the second identity} by computing
\begin{align}  \label{unitmass}
\sum_{\mu \in \mb{J}_{n}} P^{{ \psi} } _{n}(\mu) {\color{black}  \frac{ d_{\mu}^{(n)}(x_1,\ldots,x_{r_n})}{d_{\mu}^{(n)}(1,\ldots,1)} } & = \frac{ \det \left [ \sum_{k=0}^\infty \langle \tilde{P}^{(\alpha_n)}_{k} , (x-1)^{r_n-i} {{\psi} }\rangle_{\alpha_n} \tilde{P}_{k}^{(\alpha_n)}(x_j) \right ]_{i,j=1}^{r_n} }{\det [ x^{r_n -i}_j ]_{i,j=1}^{r_n}}  \\
& =\frac{ \det \left [ (x_j-1)^{r_n-i} {{ \psi} }(x_j) \right ]_{i,j=1}^{r_n} }{\det [ (x_j-1)^{r_n-i} ]_{i,j=1}^{r_n} } \notag =\prod_{j=1}^{r_n} {{ \psi} }(  x_j ) , \notag
\end{align}
where we have used a variant of the Cauchy--Binet formula (see Lemma \ref{cauchybinet}) and orthogonal decomposition \eqref{orthogonaldecomposition}. Similarly, if ${ \psi  \equiv 1}$, then
$$
\ba
P^{{ \psi}}_n(0) & = \frac{ \det \left [ \sum_{k=0}^{r_n-1} \langle \tilde{P}^{(\alpha_n)}_{k} , (x-1)^{r_n-i}  \rangle_{\alpha_n} \tilde{P}_{k}^{(\alpha_n)}(x_j) \right ]_{i,j=1}^{r_n} }{\det [ x^{r_n -i}_j ]_{i,j=1}^{r_n}}  { \Bigg|_{x_1 = \cdots = x_n=1}} \\
&=   \frac{ \det \left [ (x_j-1)^{r_n-i}  \right ]_{i,j=1}^{r_n} }{\det [ (x_j-1)^{r_n-i} ]_{i,j=1}^{r_n} }  { \Bigg|_{x_1 = \cdots = x_n=1}} = 1
\ea
$$
If instead $\lambda \neq 0$, then $P^{{ \psi}}_n(\lambda)=0$ since the first expression above will have at least one row of zeros by orthogonality to polynomials of lower degree.

Assume $n-1$ odd, so that $\alpha_{n-1} = \alpha$ and $r_n = r_{n-1}$. Then we have for $\mu \in \mb{J}_{n-1}$
\begin{equation} \label{consistencycalc}
\ba
\sum_{\lambda \in \mb{J}_n} P^\gamma_{n}(\lambda) T^n_{n-1}(\lambda, \mu) & = {\color{black} d_{\mu}^{(n)}(1,\ldots,1)} \sum_{\lambda \in \mb{J}_n} \det[ \phi_{n-1}( \tilde{\mu}_i, \tilde{\lambda}_j) ]_{i,j =1}^{r_{n}} \det \left [  \Psi^{n, \psi_\gamma}_{r_n-i}(\widetilde{\lambda}_j)  \right]_{i,j=1}^{r_n}  \\
& = {\color{black} d_{\mu}^{(n)}(1,\ldots,1)} \det \left [  (\phi_{n-1}*\Psi^{n, \psi_\gamma}_{r_n-j})(\widetilde{\mu}_i) \right ]_{i,j =1}^{r_{n}} \\
& = {\color{black} d_{\mu}^{(n)}(1,\ldots,1)} \det \left [ \Psi^{n-1, \psi_\gamma}_{r_{n-1}-j}(\widetilde{\mu}_i) \right ]_{i,j =1}^{r_{n-1}} = P^\gamma_{n-1}(\mu).
\ea
\end{equation}
where we have used Lemma \ref{cauchybinet} and the composition rule Lemma \ref{comprule}. The exact same calculations of  \eqref{consistencycalc} hold for the case $n-1$ even (where instead $\alpha_{n-1} =\alpha+1$ and $r_n = r_{n-1}-1$), except to justify its last equality, we need the fact that the $r_n$th row where the conventions $\widetilde{\mu}_{r_n} \equiv -1$ and $\phi_{n-1}(-1,\cdot) \equiv 1$ has $j$th entry
\begin{equation} \label{Muppertriangular}
(\phi_{n-1}*\Psi^{n, \psi_\gamma}_{r_n-j})(-1)  = \sum_{ s \geq 0} \Psi^{n, \psi_\gamma}_{r_n-j}(s) = \sum_{ s \geq 0} \left \langle \frac{\bar{P}^{(\alpha)}_{\widetilde{\lambda}_i}}{2^{\alpha+\beta+1}\Gamma(\alpha+1)} , (x-1)^{r_n-j} \psi_\gamma \right \rangle_{\alpha+1} \bar{\bar{P}}^{\alpha}_s(1) = \delta_{r_n,j},
\end{equation}
where we have used the fact $ \bar{\bar{P}}^{\alpha_n}_s(1) = 1$ for $n$ odd and orthogonal decomposition \eqref{orthogonaldecomposition} in the last equality. This completes the proof. 
\epf

\section{Multilevel Markov process} \label{sectionMarkov}

\subsection{Single level dynamics}

The next proposition dictates a natural way to transition between the measures $P^\gamma_n$ on a single level $\mb{J}_n$. Define $P^\psi_n$ by replacing $\psi_\gamma$ with $\psi$ in the expression for $P^\gamma_n$ of Proposition \ref{partitionmeasure}. For any fixed $\psi \in C^1[-1,1]$, define a $\mb{J}_n \times \mb{J}_n$--matrix $T^\psi_n$ with entries { (recalling the scaling \eqref{scaling})}
$$
T^\psi_{n}(\mu, \lambda) : = \det \left [ \frac{ \left \langle  \bar{\bar{P}}^{(\alpha_n)}_{\widetilde{\mu}_i} , \bar{P}^{(\alpha_n)}_{\widetilde{\lambda}_j}\psi \right \rangle_{\alpha_n}}{2^{\alpha+\beta+1}\Gamma(\alpha+1)} \right]_{i,j=1}^{r_n} \cdot \frac{{\color{black} d_{\lambda}^{(n)}(1,\ldots,1)}}{{\color{black} d_{\mu}^{(n)}(1,\ldots,1)}}.
$$ 
{For a class of functions in $C^1[-1,1]$, the matrix $T_n^{\psi}$ will be stochastic for all $n\geq 1$. In light of Theorem 5.2 of \cite{aogo1}, a natural class of functions to consider is the set $\mathcal{Y}$ of convex combinations of products of functions of the form
$$
e^{\gamma(x-1)}, \quad 1+p(x-1), \quad \frac{1}{1-q(x-1)}.
$$ 
Indeed, if $T_n^{\psi}$ is stochastic for all $n\geq 1$, then $P_n^{\psi}$ is a probability measure for all $n\geq 1$, so $\psi\in \mathcal{Y}.$ However, the converse is not true: Theorem 5.2 of \cite{aogo1} only establishes that $P_n^{\psi}$ is a probability measure for {sufficiently large} $n$, not that $T_n^{\psi}$ is stochastic {for all} $n\geq 1$. 
}

The next result generalizes statements in Section 3.1 of \cite{abjk1}.

\begin{proposition} \label{Markov1}
Let $\psi, \psi_1, \psi_2 \in C^1[-1,1]$. Let $\psi_\gamma(x) : = e^{\gamma(x-1)}$.
\begin{enumerate}

\item
$
\sum_{\mu \in \mathbb{J}_n} P^{\psi_1}_n(\mu) T^{\psi_2}_{n}(\mu,\lambda) = P^{\psi_1\cdot \psi_2}_n(\lambda),
$ 
and similarly $T^{\psi_1}_n T^{\psi_2}_n=T^{\psi_1\psi_2}_n$

\item $\sum_{\lambda \in \mb{J}_n} T^\psi_n(\mu,\lambda) = \psi(1)^{r_n}$

\item If ${ \psi_m} \in C^1[-1,1]$ converge uniformly to $\psi$, then $T^{{ \psi_m}}_{n}(\mu, \lambda) \to T^\psi_{n}(\mu, \lambda)$ as $m \to \infty$.

{\color{black} \item Set $\phi(x) = 1+a(x-1)$ where $a$ is non--negative. For $n$ odd, $T_n^{\phi}$ is stochastic if 
\small
\begin{multline*}
a^{-1} \geq \text{min}\Bigg( 2 , \frac{4 \alpha ^3+\alpha ^2 (8 \beta +9)+4 \alpha 
   \left(\beta ^2+3 \beta +1\right)+\beta  (3 \beta +4)}{(\alpha +\beta ) (\alpha +\beta +1) (\alpha +\beta +2)} ,\\
    \frac{4 \alpha ^3+8 \left(2 \alpha ^2+3 \alpha (\beta +1) 
    \right)+\alpha ^2 (8 \beta +9)+(4 \alpha +8) \left(\beta ^2+3 \beta +1\right)+24 (\alpha
   +\beta +1)+\beta  (3 \beta +4)+16}{(\alpha +\beta +2) (\alpha +\beta +3) (\alpha +\beta +4)}\Bigg),
\end{multline*}
\normalsize
In particular, this lower bound is equal to $2$ if and only if the two inequalities $(2 \alpha +2 \beta +7) \left(\alpha ^2-\beta ^2\right)\leq 0$ and $(2 \alpha +2 \beta +3) \left(\alpha ^2-\beta ^2\right)\leq 0$ hold.

For $n$ even, $T_n^{\phi}$ is stochastic if
\begin{multline*}
a^{-1} \geq \text{min}\Bigg( 2, \frac{4 \alpha ^3+\alpha ^2 (8 \beta +21)+\alpha  \left(4 \beta ^2+28
   \beta +34\right)+7 \beta ^2+24 \beta +17}{(\alpha +\beta +1) (\alpha
   +\beta +2) (\alpha +\beta +3)}, \\
\frac{4 \alpha ^3+\alpha ^2 (8 \beta +37)+\alpha  \left(4 \beta ^2+52
   \beta +114\right) +15 \beta ^2 +96\beta
   +129}{(\alpha +\beta +3) (\alpha +\beta +4) (\alpha +\beta +5)}  \Bigg),
\end{multline*}
\normalsize
In particular, the lower bound is equal to $2$ if and only if the two inequalities  $(2(\alpha+1)+3)(\alpha+1) \geq (2\beta+3)\beta$ and
$$
2 \alpha ^3+\alpha ^2 (2 \beta +13)+\alpha  \left(-2 \beta ^2+4 \beta
   +20\right)+2 \beta +9\geq \beta ^2 (2 \beta +9)
$$
hold.
}

\item For any fixed $\alpha, \beta > -1$, the semigroup $\{ T^{\psi_\gamma}_n \}_{ \gamma \geq 0}$ operating on the Banach space of absolutely summable functions $l^1(\mb{J}_n)$ on $\mb{J}_n$ is stochastic and { satisfies} $ \Vert T^{\psi_\gamma}_n - Id \Vert_{l^1(\mb{J}_n)} \overset{\gamma \to 0}{\to} 0$.

\end{enumerate}
\end{proposition}
\bpf
Both items of the first point follow readily from Lemma \ref{cauchybinet} and orthogonal decomposition \eqref{orthogonaldecomposition}; for example,
$$
\ba
\sum_{\mu \in \mathbb{J}_n} P^{\psi_1}_n(\mu) T^{\psi_2}_{n}(\mu,\lambda) &=  {\color{black} d_{\lambda}^{(n)}(1,\ldots,1)} \sum_{\mu \in \mathbb{J}_n} \det \left [  \frac{\langle \bar{P}^{(\alpha_n)}_{\widetilde{\mu}_j} , (x-1)^{r_n-i} \psi_1 \rangle_{\alpha+1}}{2^{\alpha+\beta+1}\Gamma(\alpha+1)}  \right]_{i,j=1}^{r_n} \det \left [ \frac{\left \langle  \bar{\bar{P}}_{\widetilde{\mu}_i,\alpha_n} , \bar{P}_{\widetilde{\lambda}_j,\alpha_n}\psi_2 \right \rangle_{\alpha_n} }{2^{\alpha+\beta+1}\Gamma(\alpha+1)} \right]_{i,j=1}^{r_n} \\
& =  {\color{black} d_{\lambda}^{(n)}(1,\ldots,1)} \cdot \det \left [ \left \langle \sum_{k=0}^\infty  \left \langle  \tilde{P}_{k,\alpha_n} , \frac{\bar{P}_{\widetilde{\lambda}_j,\alpha_n}\psi_2}{2^{\alpha+\beta+1}\Gamma(\alpha+1)} \right \rangle_{\alpha_n} \tilde{P}_{k,\alpha_n} , (x-1)^{r_n - i} \psi_1 \right \rangle_{\alpha_n} \right]_{i,j=1}^{r_n}  \\
& ={\color{black} d_{\lambda}^{(n)}(1,\ldots,1)} \cdot \det \left [ \frac{\left \langle \bar{P}_{\widetilde{\lambda}_j,\alpha_n}\psi_2 , (x-1)^{r_n - i} \psi_1 \right \rangle_{\alpha_n}}{2^{\alpha+\beta+1}\Gamma(\alpha+1)} \right]_{i,j=1}^{r_n} = P^{\psi_1\psi_2}_n(\lambda).
\ea
$$
The second point is similar to \eqref{unitmass}: 
\begin{align*}
\sum_{\lambda \in \mb{J}_{n}} T^\psi_n(\mu, \lambda) & = \frac{1  }{{\color{black} d_{\mu}^{(n)}(1,\ldots,1)}} \sum_{\lambda \in \mb{J}_{n}} \det \left [ \frac{ \left \langle  \bar{\bar{P}}^{(\alpha_n)}_{\widetilde{\mu}_i} , \bar{P}^{(\alpha_n)}_{\widetilde{\lambda}_j}\psi \right \rangle_{\alpha_n}}{2^{\alpha+\beta+1}\Gamma(\alpha+1)} \right]_{i,j=1}^{r_n}  \frac{\det \left [\bar{\bar{P}}_{\widetilde{\lambda}_i}^{(\alpha_n)}(x_j) \right ]_{i,j=1}^{r_n}}{\det [ x_j^{r_n-i} ]_{i,j=1}^{r_n}}  { \Bigg|_{x_1 = \cdots = x_n=1}}   \\
& =\frac{1}{{\color{black} d_{\mu}^{(n)}(1,\ldots,1)}} \frac{\det \left [ \sum_{k=0}^\infty \langle \bar{\bar{P}}_{\widetilde{\mu}_i}^{(\alpha_n)} \psi,  \tilde{P}_{k}^{(\alpha_n)} \rangle_{\alpha_n}\tilde{P}_{k}^{(\alpha_n)}(x_j) \right]_{i,j=1}^{r_{n}}}{\det [ x_j^{r_n-i} ]_{i,j=1}^{r_n}}  { \Bigg|_{x_1 = \cdots = x_n=1}} \\
& =\frac{1 }{{\color{black} d_{\mu}^{(n)}(1,\ldots,1)}} \frac{\det \left [ \bar{\bar{P}}_{\widetilde{\mu}_i,\alpha_n}(x_j) \right]_{i,j=1}^{r_{n}}}{\det [ x_j^{r_n-i} ]_{i,j=1}^{r_n}}\prod_{k=1}^{r_n} \psi(x_k)  { \Bigg|_{x_1 = \cdots = x_n=1}} = \psi(1)^{r_n}
 \end{align*}
where Lemma \ref{cauchybinet} and orthogonal decomposition \eqref{orthogonaldecomposition} were used in the second and third equalities. The third point of the proposition is a straightforward application of the dominated convergence theorem. 

For the { fourth} point, { it suffices by the second point to show that the entries are non--negative}. Recall that we write $\bar{P}^{(\alpha_n)}_k = \bar{c}_k^n P^{(\alpha_n)}_k$ and $\bar{\bar{P}}^{(\alpha_n)}_k =\bar{\bar{c}}_k^n P^{(\alpha_n)}_k$. Then using the three term recurrence \eqref{recurrence} and orthogonality relations, we then compute
\begin{equation} \label{recurrencejumprate}
\ba
\frac{\langle \bar{\bar{P}}^{(\alpha_n)}_{\mu} , \bar{P}^{(\alpha_n)}_{\lambda} { \phi} \rangle_{\alpha_n}}{2^{\alpha+\beta+1} \Gamma(\alpha+1)}  = (1-{ a}) \delta_{\mu, \lambda}  &+ { a} \frac{\bar{c}_{\lambda}^n}{\bar{c}_{\mu}^n} [ A^{{ (\alpha_n, \beta)}}_{\lambda} \delta_{\mu, \lambda} +  B^{{ (\alpha_n, \beta)}}_{\lambda} \delta_{\mu, \lambda-1} +  C^{{ (\alpha_n, \beta)}}_{\lambda} \delta_{\mu, \lambda+1} ]
\ea
\end{equation}
If $\widetilde{\lambda}_i > \widetilde{\mu}_i+1$ for some $1 \leq i \leq r_n$, then $\widetilde{\lambda}_k > \widetilde{\mu}_l+1$ for $1 \leq k \leq i \leq l \leq r_n$, which implies $T^{{ \phi}}_{n}(\mu,\lambda) = 0$ (the resulting matrix admits a $2 \times 2$ block form with an off-diagonal block of $0$'s and a diagonal block with a zero vector). The same conclusion of course holds if $\widetilde{\mu}_i>\widetilde{\lambda}_i +1$, so assume $|\widetilde{\lambda}_i - \widetilde{\mu}_i|\leq 1$ for all $1 \leq i \leq r_n$. If $|\widetilde{\mu}_i - \widetilde{\mu}_{i+1}|>1$ for some $1 \leq i < r_n$, then $|\widetilde{\mu}_i - \widetilde{\lambda}_{i+1}| \leq 1$ implies $|\widetilde{\lambda}_{i}-\widetilde{\mu}_{i+1}| > 1$ and similarly $|\widetilde{\lambda}_{i}-\widetilde{\mu}_{i+1}| \leq 1$ implies $|\widetilde{\mu}_i - \widetilde{\lambda}_{i+1}| > 1$. In either case, $T^{{ \phi}}_{n}(\mu,\lambda)$ breaks into a product of determinants, one of which is the $i \times i$ upper left corner minor of the matrix defining $T^{{ \phi}}_{n}(\mu, \lambda)$. Iterating this argument reduces consideration to the case where $|\widetilde{\mu}_i - \widetilde{\mu}_{i+1}|\leq1$ for all $1 \leq i < r_n$, which means $\mu_i$ are all equal to some $p \in \Z_{\geq 0}$. Now the two blocks corresponding to $\{i : \lambda_i=p\pm1\}$ are triangular with nonnegative entries (since $B^{(\alpha_n,\beta)}_k, C^{(\alpha_n,\beta)}_k \geq 0$)  and straddle a tridiagonal block corresponding to $\{ i:\lambda_i=p = \mu_i \}$. Write $q^*:= \max\{ i:\lambda_i=p = \mu_i \}$, $q_*:= \min\{ i:\lambda_i=p = \mu_i \}$, $q:=q^*-q_*+1$, and set $\mu:= \tilde{\mu}_{q_*}$. The last argument has further reduced consideration to the determinant of the $q \times q$ tridiagonal matrix
\begin{equation} \label{tridiagonal1}
\begin{bmatrix}
1-{ a} + { a} A^{(\alpha_n,\beta)}_{\mu} &  { a} \frac{ \bar{c}_{\mu-1}^n}{\bar{c}_{\mu}^n}C^{(\alpha_n,\beta)}_{\mu-1} & 0 & 0  \\
{ a} \frac{ \bar{c}_{\mu}^n}{\bar{c}_{\mu-1}^n}B^{(\alpha_n,\beta)}_{\mu} & 1-{ a} + { a} A^{(\alpha_n,\beta)}_{\mu-1}  &  { a}  \frac{\bar{c}_{\mu-2}^n}{\bar{c}_{\mu-1}^n}C^{(\alpha_n,\beta)}_{\mu-2}  & 0    \\
0 & { a} \frac{ \bar{c}_{\mu-1}^n}{\bar{c}_{\mu-2}^n}B^{(\alpha_n,\beta)}_{\mu-1} &  1-{ a} + { a} A^{(\alpha_n,\beta)}_{\mu-2} &  \ddots   \\
\vdots & \vdots  & \ddots  & \ddots    \\
\end{bmatrix}.
\end{equation}
{\color{black} In general, if a tridiagonal matrix has the property that each diagonal entry of a given row is greater than or equal to the sum of the off-diagonals in that row, then all of its principal minors are non--negative (see page 5 of \cite{FJ}). In particular, its determinant is non--negative.} This leads us to show that the expressions
\begin{equation} \label{nonnegativitycondition}
{ a} \left [  \frac{ \bar{c}_{1}^n}{\bar{c}_{0}^n}B^{(\alpha_n,\beta)}_{1} - A^{(\alpha_n,\beta)}_0 { + 1} \right ], \ \ \ { a} \left [ \frac{ \bar{c}_{k+1}^n}{\bar{c}_{k}^n}B^{(\alpha_n,\beta)}_{k+1}+ \frac{ \bar{c}_{k-1}^n}{\bar{c}_{k}^n}C^{(\alpha_n,\beta)}_{k-1} - A^{(\alpha_n,\beta)}_{k} { + 1}  \right ], \ \ \ k \geq 1,
\end{equation}
are {\color{black} less than or equal to $1$. Since the formula defining $C^{(\alpha_n,\beta)}_{k-1}$ is equal to zero when $k=0$, it suffices to focus} on the second expression. {\color{black} When $n$} is odd, {\color{black} we} compute 
$$
\ba
&  \frac{ \bar{c}_{k+1}^n}{\bar{c}_{k}^n}B^{(\alpha_n,\beta)}_{k+1}+ \frac{ \bar{c}_{k-1}^n}{\bar{c}_{k}^n}C^{(\alpha_n,\beta)}_{k-1} - A^{(\alpha_n,\beta)}_{k} = \\
& \frac{k+\alpha+\beta+1}{2k+\alpha+\beta+1}  \frac{2(k+\alpha+1)}{2k+\alpha+\beta+2} + \frac{k+\beta}{2k+\alpha+\beta}  \frac{2k}{2k+\alpha+\beta+1} - \frac{\beta^2 - \alpha^2}{(2k + \beta +  \alpha)(2k+\alpha+\beta+1)}
 \ea
$$
{\color{black} With some algebra, it can be shown that the desired inequality is equivalent to the condition that $a^{-1}$ is bounded below by  
\footnotesize
$$
\frac{4 \alpha ^3+\alpha ^2 (8 \beta +9)+(4 \alpha +8k)
   \left(\beta ^2+3 \beta +1\right)+\beta  (3 \beta +4)+16
   k^3+24 k^2 (\alpha +\beta +1)+8 k \left(2 \alpha ^2+3 \alpha 
   (\beta +1)\right)}{(\alpha +\beta +2 k) (\alpha +\beta +2 k+1) (\alpha +\beta
   +2 k+2)}.
$$
\normalsize
By taking a derivative with respect to $k$, it can be seen that this quantity as a function of $k$ can only have real critical points if $k$ is a real solution of a cubic polynomial, which can be solved explicitly to find
$$
k=\frac{1}{8} \left(-4 \alpha -4 \beta +\sqrt[3]{7+4 \sqrt{3}}+\frac{1}{\sqrt[3]{7+4 \sqrt{3}}}-5\right) \approx \frac{1}{8}(-2.17836-4\alpha-4\beta),
$$ 
which is less than $1$ for $\alpha_n,\beta > -1$.
Therefore, this quantity is monotonic on $k\in [1,\infty)$, so it suffices to evaluate it at $k=\infty,0,1,$ which yields the three terms in the lower bound in the fourth statement. Setting the two terms from $k=0,1$ equal to $2$ yields the two inequalities $(2 \alpha +2 \beta +7) \left(\alpha ^2-\beta ^2\right)\leq 0$ and $(2 \alpha +2 \beta +3) \left(\alpha ^2-\beta ^2\right)\leq 0$.

When $n$ is even, 
\footnotesize
$$
\ba
&  \frac{ \bar{c}_{k+1}^n}{\bar{c}_{k}^n}B_{k+1}^{(\alpha,\beta)}+ \frac{ \bar{c}_{k-1}^n}{\bar{c}_{k}^n}C_{k-1}^{(\alpha,\beta)} - A_{k}^{(\alpha,\beta)} = \\
& \frac{2(k+1)(k+\beta+1)}{ (2k+\alpha+\beta+3)(2k+\alpha+\beta+2)} + \frac{2(\alpha+k+1)(k+\alpha+\beta+1)}{(2k+\alpha+\beta+2)(2k+\alpha+\beta+1)}- \frac{\beta^2 - \alpha_n^2}{(2k + \beta +  \alpha_n)(2k+\alpha_n+\beta+1)}.
 \ea
$$
\normalsize
With some algebra, $a^{-1}$ is bounded below by
\footnotesize
$$
\frac{4 \alpha ^3+\alpha ^2 (8 \beta +21)+\alpha  \left(4 \beta ^2+28
   \beta +34\right)+7 \beta ^2+24 \beta +16 k^3+24 k^2 (\alpha +\beta
   +2)+8 k \left(2 \alpha ^2+\alpha  (3 \beta +7)+\beta ^2+6 \beta
   +6\right)+17}{(\alpha +\beta +2 k+1) (\alpha +\beta +2 k+2) (\alpha
   +\beta +2 k+3)}
$$
\normalsize
Again, this quantity as a function of $k$ can only have real critical points at 
$$
k=\frac{1}{8} \left(-4 \alpha -4 \beta +\sqrt[3]{7+4
   \sqrt{3}}+\frac{1}{\sqrt[3]{7+4 \sqrt{3}}}-9\right),
$$
which is again less than $1$ for $\alpha,\beta>-1$. Plugging in $k=0,1$ yields the two quantities in the lower bound, and setting those quantities equal to $2$ yields the two inequalities.

}

To conclude the {last statement}, { note that by the third and first points 
\begin{equation} \label{limittransition}
T^{\psi_\gamma}_n =  \lim_{m \to \infty} T^{(1+\gamma(x-1)/m)^m}_n = \lim_{m \to \infty} T^{1+\gamma(x-1)/m}_n \cdots T^{1+\gamma(x-1)/m}_n. 
\end{equation}
By the second point and fourth points respectively, the entries of $T_n^{\psi_{\gamma}}$ have rows which sum to $1$ and are positive. Additionally,
\begin{equation} \label{equivform}
\Vert T^{\psi_\gamma}_n - Id \Vert_{l^1(\mb{J}_n)} = \sup_{\lambda \in \mb{J}_n} \sum_{\mu \in \mb{J}_n} (\delta_{\lambda \mu} - T^{\psi_\gamma}_n(\lambda,\mu)) = \sup_{\lambda \in \mb{J}_n} (2 - 2T_n^{\psi_\gamma}(\lambda,\lambda)).
\end{equation}
}
Note by \eqref{limittransition} that $T^{\psi_\gamma}_n(\lambda,\lambda) \geq \lim_{m \to \infty} [T^{(1+\gamma(x-1)/m)}_n(\lambda,\lambda)]^m$. The identity \eqref{equivform} then tells us it is sufficient to show $\inf_{\lambda \in \mb{J}_n} T^{1-\gamma/m+x\cdot\gamma/m}_n(\lambda,\lambda)$ goes to $1$ as $\gamma \to 0$. But our argument above tells us that $T^{1-\gamma/m+x\cdot\gamma/m}_n(\lambda,\lambda)$ breaks apart into a product of determinants of $q \times q$ matrices of the form \eqref{tridiagonal1}. Thus we are interested in
$$
\inf_{\mu \geq q} \det \left [ \left (1-{ a} \right) I_q + 
  { a^q} \cdot \begin{bmatrix}
  A_{\mu} &   \frac{ \bar{c}_{\mu-1}^n}{\bar{c}_{\mu}^n}C_{\mu-1} & 0 & 0  \\
 \frac{ \bar{c}_{\mu}^n}{\bar{c}_{\mu-1}^n}B_{\mu} &  A_{\mu-1}  &    \frac{\bar{c}_{\mu-2}^n}{\bar{c}_{\mu-1}^n}C_{\mu-2}  & 0    \\
0 &  \frac{ \bar{c}_{\mu-1}^n}{\bar{c}_{\mu-2}^n}B_{\mu-1} &   A_{\mu-2} &  \ddots   \\
\vdots & \vdots  & \ddots  & \ddots    \\
\end{bmatrix} \right ].
$$
For large enough $m$, this infimum is nonzero and in fact positive by our reasoning above. By continuity, it goes to $1$ as $\gamma \to { 0}$, as required.

\epf

{\color{black} 
\begin{remark}
In Proposition 3.3 of \cite{abjk1}, the values $\alpha = -1/2$ and $\beta=-1/2$ are considered with the conclusion that $T_n^{\phi}$ is stochastic for $\phi(x) = 1+a(x-1)$ when $a$ is bounded above by $1/2$. Therefore, Proposition \ref{Markov1}.4 is a generalization of that result.

In general, there will be values of $a$ for which $T^{\phi}_n$ with $\phi(x)=1+a(x-1)$ will not be stochastic. For example, consider $T_n^{\phi}$ for $n=1,2$. By \eqref{recurrence},
$$
(1+a(x-1)) P_k^{(\alpha)}(x) = (1+a(A_k^{(\alpha,\beta)}-1))P_k^{(\alpha)}(x) + aB_k^{(\alpha,\beta)}P_{k-1}^{(\alpha)}(x) + aC_k^{(\alpha,\beta)}P_{k+1}^{(\alpha)}(x).
$$
Since $B_k^{(\alpha,\beta)},C_k^{(\alpha,\beta)}$ are positive and $ 1 > \tfrac{\vert \beta-\alpha\vert}{\alpha+\beta+2} = \vert A_0\vert \geq \vert A_1 \vert \geq \ldots$ with $\lim_{k\rightarrow\infty} \vert A_k\vert =0$, then $T_n^{\phi}$ will be stochastic if and only if
$$
a \leq \frac{1}{1-A_0^{(\alpha_n,\beta)}} = \frac{\alpha_n+\beta+2}{2\alpha_n+2 }.
$$

Due to results of \cite{abic}, it might be possible to interpret $T_n^{\phi}$ as a discrete--time particle system with Bernoulli jumps, but this direction is not pursued here.
\end{remark}
}
\begin{remark}
When $\psi=(1-q(x-1))^{-1}$, the matrix $T_n^{\psi}$ has been previously considered in \cite{jk1,jk2} with the values $\vert \alpha \vert = \vert \beta \vert = 1/2$ arising from the representation theory of the orthogonal and {symplectic} groups. In those cases, there was the property that $T_n^{\psi}(\lambda,\mu)$ was nonzero only if there existed a $\nu$ such that $\nu\prec \lambda,\mu$, which lead to the corresponding Markov process being interpreted as a discrete--time particle system having geometric jump rates with states $\nu$ at half--integer times. However, one can check computationally that this property will not hold for general $\alpha,\beta>-1$.
\end{remark}

\subsection{Intertwined multilevel dynamics: discrete steps}

We now implement the procedure for constructing multilevel dynamics with certain restrictions introduced in Borodin--Ferrari \cite{abpf1}, based on the construction in Diaconis--Fill \cite{pdjf1}. These stochastic operators $T^n_{n-1}$ defined above ensure that the interlacing condition is preserved when we link single level dynamics to a multilevel evolution, and the following intertwining relationship is key.

\begin{proposition} \label{interlacingrelations}
Fix $\psi \in C^1([-1,1])$ with $\psi(1)=1$. For $n \geq 1$, the stochastic operators $T^\psi_n$ and $T_n^{n+1}$ satisfy the intertwining relations 
$$
\Delta_n^{n+1} : =T^{n+1}_n T^\psi_n = T^\psi_{n+1}T^{n+1}_n.
$$

\end{proposition}
\bpf
First assume $n$ is even, so that $\alpha_n = \alpha+1$, $\alpha_{n+1} = \alpha$, $r_n = r_{n+1} -1$, and $\phi_n(s,t) = \phi_n(s) 1_{(s<t)}$. For $\mu \in \mb{J}_{n+1}$ and $\lambda \in \mb{J}_n$, Lemma \ref{indication} allows us to write
$$
\ba
(T^{n+1}_nT^\psi_n)(\mu, \lambda) & = \frac{{\color{black} d_{\lambda}^{(n)}(1,\ldots,1)}}{{\color{black} d_{\mu}^{(n+1)}(1,\ldots,1)}} \sum_{z \in \mb{J}_n} \det \left [ \phi_n(\widetilde{z}_j,\widetilde{\mu}_i) \right ]_{i,j=1}^{r_{n+1}} \det \left [ \frac{ \left \langle  \bar{\bar{P}}_{\widetilde{z}_j}^{(\alpha+1)}, \bar{P}_{\widetilde{\lambda}_j}^{(\alpha+1)} \psi \right \rangle_{\alpha+1}}{2^{\alpha+\beta+1} \Gamma(\alpha+1)} \right ]_{i,j=1}^{r_{n}}.
\ea
$$
Expand the first determinant along the $r_{n+1}$-column, where the convention $z_{r_{n+1}} \equiv -1$ applies. For each $1 \leq l \leq r_{n+1}$, the $l$th resulting summand is
\begin{align*}
 (-1)^{r_{n+1}+l} \sum_{z \in \mb{J}_n} \det & \left [ \phi_n(\widetilde{z}_j, \widetilde{\mu}_i) \right ]_{\substack{1 \leq i \neq l \leq r_{n+1} \\ 1 \leq j \leq r_n}}  \det \left [ \frac{1}{2^{\alpha+\beta+1} \Gamma(\alpha+1)} \left \langle \tilde{P}_{\widetilde{z}_j}^{(\alpha+1)}, \tilde{P}_{\widetilde{\lambda}_j}^{(\alpha+1)} \psi \right \rangle_{\alpha+1} \right ]_{i,j=1}^{r_{n}} \\
& = (-1)^{r_{n+1}+l} \det \left [ \frac{1}{2^{\alpha+\beta+1} \Gamma(\alpha+1)} \left \langle \sum_{k = 0}^{\widetilde{\mu}_i - 1} \bar{P}_{k}^{(\alpha)}, \bar{P}_{\widetilde{\lambda}_j}^{(\alpha+1)} \psi \right \rangle_{\alpha+1} \right ]_{\substack{1 \leq i \neq l \leq r_{n+1} \\ 1 \leq j \leq r_n}}  \\
& =(-1)^{r_{n+1}+l}  \det \left [ \frac{1}{2^{\alpha+\beta+1} \Gamma(\alpha+1)} \left \langle \frac{\bar{\bar{P}}_{\widetilde{\mu}_i}^{(\alpha)}-1}{x-1}, \bar{P}_{\widetilde{\lambda}_j}^{(\alpha+1)} \psi \right \rangle_{\alpha+1} \right ]_{\substack{1 \leq i \neq l \leq r_{n+1} \\ 1 \leq j \leq r_n}}  \\
& =  (-1)^{r_{n+1}+l} \det \left [ \frac{1}{2^{\alpha+\beta+1} \Gamma(\alpha+1)} \left \langle (\bar{\bar{P}}_{\widetilde{\mu}_i}^{(\alpha)}\psi)(1) - \bar{\bar{P}}_{\widetilde{\mu}_i}^{(\alpha)}\psi, \phi_n(\widetilde{\lambda}_j) \bar{\bar{P}}_{\widetilde{\lambda}_j}^{(\alpha+1)} \right \rangle_{\alpha} \right ]_{\substack{1 \leq i \neq l \leq r_{n+1} \\ 1 \leq j \leq r_n}}  \\
 & =  (-1)^{r_{n+1}+l} \det \left [ \frac{1}{2^{\alpha+\beta+1} \Gamma(\alpha+1)} \left \langle \bar{\bar{P}}_{\widetilde{\mu}_i}^{(\alpha)} \psi, \phi_n(\widetilde{\lambda}_j) \sum_{k = \widetilde{\lambda}_j+1}^\infty    \bar{P}_{k}^{( \alpha)} \right \rangle_{\alpha} \right ]_{\substack{1 \leq i \neq l \leq r_{n+1} \\ 1 \leq j \leq r_n}} 
 \\
 & =(-1)^{r_{n+1}+l}  \det \left [ \frac{1}{2^{\alpha+\beta+1} \Gamma(\alpha+1)} \left \langle \bar{\bar{P}}_{\widetilde{\mu}_i}^{( \alpha)}\psi, \sum_{k =0}^\infty \phi_n(\widetilde{\lambda}_j, k) \bar{P}_{k}^{(\alpha)} \right \rangle_{\alpha} \right ]_{\substack{1 \leq i \neq l \leq r_{n+1} \\ 1 \leq j \leq r_n}},
\end{align*}
where the first equality follows from Proposition \ref{identities}, the second from Lemma \ref{cauchybinet}, and the third from the fact $(\bar{\bar{P}}_{\widetilde{\mu}_i}^{(\alpha)}\psi)(1) =1$ along with basic properties of determinants. Summing this calculation over $1 \leq l \leq r_n$ gives
\begin{align*}
(T^{n+1}_nT^\psi_n)(\mu, \lambda) & = \frac{{\color{black} d_{\lambda}^{(n)}(1,\ldots,1)}}{{\color{black} d_{\mu}^{(n+1)}(1,\ldots,1)}}   \sum_{l=1}^{r_{n+1}} (-1)^{r_{n+1}+l} \det \left [  \frac{1}{2^{\alpha+\beta+1} \Gamma(\alpha+1)} \left \langle \bar{\bar{P}}_{\widetilde{\mu}_i}^{(\alpha)}\psi, \sum_{k =0}^\infty \phi_n(\widetilde{\lambda}_j, k) \bar{P}_{k}^{(\alpha)} \right \rangle_{\alpha} \right ]_{\substack{1 \leq i \neq l \leq r_{n+1} \\ 1 \leq j \leq r_n}} \\
& = \frac{{\color{black} d_{\lambda}^{(n)}(1,\ldots,1)}}{{\color{black} d_{\mu}^{(n+1)}(1,\ldots,1)}}  \det \left [  \frac{1}{2^{\alpha+\beta+1} \Gamma(\alpha+1)} \left \langle \bar{\bar{P}}_{\widetilde{\mu}_i}^{(\alpha)} \psi, \sum_{k =0}^\infty \phi_n(\widetilde{\lambda}_j, k) \bar{P}_{k}^{( \alpha)} \right \rangle_{\alpha} \right ]_{i,j=1}^{r_{n+1}} \\
& = \frac{{\color{black} d_{\lambda}^{(n)}(1,\ldots,1)}}{{\color{black} d_{\mu}^{(n+1)}(1,\ldots,1)}}  \sum_{z \in \mb{J}_{n+1}} \det \left [  \frac{1}{2^{\alpha+\beta+1} \Gamma(\alpha+1)} \left \langle \bar{\bar{P}}_{\widetilde{\mu}_i}^{(\alpha)} \psi,  \bar{P}_{z_j}^{(\alpha)} \right \rangle_{\alpha} \right ]_{i,j=1}^{r_{n+1}} \cdot \det  \left [ \phi_n(\widetilde{\lambda}_j, z_i)  \right ]_{i,j=1}^{r_{n+1}} \\
& = (T^\psi_{n+1} T^{n+1}_n)(\mu,\lambda),
\end{align*}
where the second equality follows because the entries with $j=r_{n+1}$ are all $1$ because of the combination of the convention $\widetilde{\lambda}_{r_{n+1}} \equiv -1$ and the fact $(\bar{\bar{P}}_{\widetilde{\mu}_i}^{( \alpha)}\psi)(1) =1$ along with orthogonal decomposition \eqref{orthogonaldecomposition}. 

The much more straightforward ``$n$ odd'' case involves the other components of Proposition \ref{identities} and \ref{indication}, and is left to the reader.

\epf

\subsection{Description of continuous-time multilevel dynamics} \label{continuoustime}
For $u,t \in  \mathbb{J}_{n,paths}$ define
\begin{equation} \label{conditionalprob}
L^{k+1, \psi}_k(u,t):=\frac{T^\psi_{k+1}(u^{k+1}, t^{k+1}) T_k^{k+1}(t^{k+1}, t^k)}{\Delta^{k+1}_k(u^{k+1},t^k)} 1_{(\Delta^{k+1}_k(u^{k+1},t^k) \neq 0 )}.
\end{equation}
In words, it is the probability of the transition ``$u^{k+1} \to t^{k+1} \to t^k$" by a jump and cotransition conditional on ``$u^{k+1} \to t^k$" occurring by such steps; note this quantity only depends on $u^{k+1},t^k,t^{k+1}$. Just as the stochastic (for suitable $\psi$) operators $T_n^\psi$ account for transitions between the probability measures $P_n^\psi$ on $\mb{J}_{n}$, the stochastic (by Proposition \ref{Markov1}) transition operator
$$
A_n^{\psi}(u, t) :=  T^\psi_1(u^1,t^1) \cdot \prod_{k=1}^{n-1} L^{k+1, \psi}_k(u,t),
$$
supplies an evolution of signed measures on $\mathbb{J}_{n,paths}$ of the form
\begin{equation} \label{seqmeasure}
P_{n,paths}^{\psi}(t): = P^{\psi}_n(t^n) \cdot \prod_{k=1}^{n-1} T^{k+1}_k(t^{k+1},t^k), \ \ \ t \in \mb{J}_{n,paths}
\end{equation}
in the sense 
\begin{equation} \label{pathtrans}
(P_{n,paths}^{\psi_1} A^{\psi_2}_n)(t) = P_{n,paths}^{\psi_1\psi_2}(t), \ \  \psi_1, \psi_2 \in C^1[-1,1], \ \ \ t \in \mb{J}_{n,paths}.
\end{equation}
This identity can be proven by induction, using the fact
$$
\sum_{\mb{J}_k \ni u^{k} \prec u^{k+1} \in \mb{J}_{k+1}} T^{k+1}_k(u^{k+1},u^k) L^{k+1, \psi_2}_k(u,t)  = T^{k+1}_k(t^{k+1},t^k), 
$$
and, by the first part of Proposition \ref{Markov1}, the fact
$$
\sum_{u^n \in \mb{J}_n} P^{\psi_1}_n(u^n)  T^{\psi_2}_k(u^{n},t^{n} ) = P^{\psi_1\psi_2}_n(t^n). 
$$

Now let $\mf{B}_n$ be the Banach space defined as the completion of the subspace of $l^1(\mb{J}_{n,paths})$ consisting of measures of the form \eqref{seqmeasure} corresponding to functions $\psi$ {\color{black} $\in C^1[-1,1]$}. The stochastic operators $\{ A^{\psi_\gamma}_n \}_{\gamma \geq 0}$ form a semigroup on $\mf{B}_n$ by \eqref{pathtrans}. The {uniform continuity property} for $\{ A^{\psi_\gamma}_n \}_{\gamma \geq 0}$ follows from the same property of $\{ T^{\psi_\gamma}_n \}_{\gamma \geq 0}$ (the fifth part of Proposition \ref{Markov1}): the form  \eqref{seqmeasure} implies
$$
\ba
\sum_{u \in \mb{J}_{n,paths}} |P_{n,paths}^{\psi_1 \psi_2}(u) - P_{n,paths}^{\psi_1}(u)| & = \sum_{u^n \in \mb{J}_{n}} |P_{n}^{\psi_1 \psi_2}(u^n) - P_{n}^{\psi_1}(u^n)| \sum_{\substack{u \in \mb{J}_{n-1,paths} \\ u^{n-1} \prec u^n}} \prod_{k=1}^{n-1} T^{k+1}_k(u^{k+1},u^k) \\
& = \sum_{u^n \in \mb{J}_{n}} |P_{n}^{\psi_1 \psi_2}(u^n) - P_{n}^{\psi_1}(u^n)|,
\ea
$$
so we have  
$$
\Vert A^{\psi_\gamma}_n - Id \Vert_{\mf{B}_n} = \Vert T^{\psi_\gamma}_n - Id \Vert_{l^1(\mb{J}_n)} \overset{\gamma \to 0}{\to} 0.
$$
Hence { (see, e.g. Theorem 2.6.1 of \cite{HP} or Proposition 3.13 of \cite{abjk1} )  }, there exists an operator $\bar{Q}_n$ on $\mf{B}_n$ such that $P_{n,paths}^{\psi} e^{\gamma \bar{Q}_n} = P_{n,paths}^{\psi} A^{\psi_\gamma}_n$, so we endeavor to arrive at a description of $\bar{Q}_n = \frac{d}{d\gamma} \big |_{\gamma=0} A^{\psi_\gamma}_n$. This involves computing
\begin{equation} \label{component}
\frac{T^{\psi_\gamma}_k(u^{k},t^{k} ) T_{k-1}^{k}(t^{k}, t^{k-1})}{\Delta^{k}_{k-1}(u^{k},t^{k-1})} 
\end{equation}
to second order (note in this expression the ratios of the ``$d_n(\cdot)$" terms and the factors ``$ \prod_{\ell=1}^{r_{k-1}} \phi_{k-1}(\tilde{t}_\ell^{k-1})$" cancel). By construction, we arrive at the same particle pattern and push--block dynamics as in the orthogonal and {symplectic} cases, but with different jump rates. More precisely, if a particle is at position $\lambda \geq 0$ on level $n \geq 1$ and is allowed to jump to position $\mu$, then it will do so at a rate
$$
 \frac{ \left \langle  \bar{\bar{P}}^{(\alpha_n)}_{\lambda} , \bar{P}^{(\alpha_n)}_{\mu}x \right \rangle_{\alpha_n}}{2^{\alpha+\beta+1}\Gamma(\alpha+1)} = 
\begin{cases}
\frac{\bar{c}_{\lambda+1}}{\bar{c}_{\lambda}}  B^{(\alpha_n,\beta)}_{\lambda+1} = \frac{\bar{\bar{c}}_{\lambda}^n}{\bar{\bar{c}}_{\lambda+1}^n} C^{(\alpha_n,\beta)}_{\lambda}, & \lambda \geq 0, \mu = \lambda+1 \\
\frac{\bar{c}_{\lambda-1}}{\bar{c}_{\lambda}} C^{(\alpha_n,\beta)}_{\lambda-1} = \frac{\bar{\bar{c}}^n_{\lambda}}{\bar{\bar{c}}^n_{\lambda-1}}  B^{(\alpha_n,\beta)}_{\lambda}, & \lambda \geq 1, \mu = \lambda -1 \\
\end{cases}.
$$

Let us exemplify the calculation. Assume the system is in a state so that a wall jump is possible at an odd level $k \geq 1$, i.e., a state $t \in \mb{J}_{n,paths}$ that satisfies $t^k_{r_k} = \widetilde{t^k_{r_k}} = 0$, $t^k_{r_{k-1}} \geq 1$, so $\widetilde{t^k_{r_{k-1}}} \geq 2$, and if $k \geq 3$, $t^{k-1}_{r_{k-1}} = \widetilde{t^{k-1}_{r_{k-1}} } \geq 1$. Consider a transition to $u \in \mb{J}_{n,paths}$ that agrees with $t \in \mb{J}_{n,paths}$ except that $u^k_{r_k} = \widetilde{u^k_{r_k}}= 1$. Noting that $\psi_\gamma(x)=1 + \gamma(x-1)+O(\gamma^2)$ for small $\gamma >0$, we use \eqref{recurrencejumprate} to compute, up to second-order in $\gamma$, (omitting the superscript ``$k$")
$$
\ba
& \det  \left [ \frac{1}{2^{\alpha+\beta+1} \Gamma(\alpha+1)} \left \langle \bar{P}_{\widetilde{t}_i}^{(\alpha)} , \bar{\bar{P}}_{\widetilde{u}_j}^{(\alpha)} \psi_\gamma \right \rangle_{\alpha} \right ]_{i,j = 1}^{r_k} \\
& \approx \det \begin{bmatrix}
1-\gamma + \gamma A_{\widetilde{t_1}} &  \gamma \frac{ \bar{c}_{\tilde{t}_1-1}^n}{\bar{c}_{\tilde{t}_1}^n}C_{\tilde{t}_1-1} \delta_{\widetilde{t_1}, \widetilde{t_2} + 1} & 0 & 0  &0 \\
\gamma \frac{ \bar{c}_{\widetilde{t_2}+1}^n}{\bar{c}_{\widetilde{t_2}}^n}B_{\widetilde{t_2}+1} \delta_{\widetilde{t_2}+1, \widetilde{t_1}} & 1-\gamma + \gamma A_{\widetilde{t_2}}  &  \gamma  \frac{\bar{c}_{\widetilde{t_2}-1}^n}{\bar{c}_{\widetilde{t_2}}^n}C_{\widetilde{t_2}-1} \delta_{\widetilde{t_2}, \widetilde{t_3} + 1}  & 0 & 0    \\
\vdots & \vdots  & \ddots   & \vdots & 0  \\
0 & \cdots  & \cdots   &  (1-\gamma) + \gamma A_{\widetilde{t}_{r_{k-1}}} &  \gamma \frac{\bar{c}^n_{\tilde{t}_{r_{k-1}-1}}}{\bar{c}^n_{\tilde{t}_{r_{k-1}}}} C_1 \delta_{ \tilde{t}_{r_{k-1}}, 2}  \\
0 & \cdots  & \cdots   & 0 & \gamma \frac{\bar{c}_{1}}{\bar{c}_{0}}  B^{(\alpha_n,\beta)}_{1}  \\
\end{bmatrix}  \\
& = \gamma \cdot \frac{2(\alpha+1)}{(\alpha+\beta+2)} + O(\gamma^2).
\ea
$$
Similar calculations in the denominator of \eqref{component} yield an order of $1 + O(\gamma)$. Since the terms ``$d_{\lambda}^{(n)}(\cdot)$" occurring in the numerator and denominator all cancel, we arrive at
$$
\ba
\frac{T^{\phi_\gamma}_k(t^{k},u^{k} ) T_{k-1}^{k}(u^{k}, u^{k-1})}{\Delta^{k}_{k-1}(u^{k},t^{k-1})} &= \frac{\gamma \cdot \frac{2(\alpha+1)}{(\alpha+\beta+2)} + O(\gamma^2)}{ 1 + O(\gamma)} \\
& = \left [ \gamma \cdot \frac{2(\alpha+1)}{(\alpha+\beta+2)} + O(\gamma^2)\right ] (1 + O(\gamma)) = \gamma \cdot \frac{2(\alpha+1)}{(\alpha+\beta+2)} + O(\gamma^2).
\ea
$$
where the second equality uses the geometric series. The other cases can be handled similarly.

\section{Determinantal correlations} \label{sectiondeterminantal}

\subsection{Probability measures on partition paths} \label{sectioncentral}

For $u \in \mb{J}_{n,seq}$, consider the cylinder sets
$
C_u: = \{ t \in \mb{J}_{\infty, seq} | t^1 = u^1, \ldots, t^n = u^n \}.
$
We now realize the measures $P^\gamma_n$ on $\mathbb{J}_n$ of Proposition \ref{partitionmeasure} as embedded in a single probability measure $P^\gamma$  on $\mb{J}_{\infty, seq}$ by the prescription
$$
P^\gamma(C_u): = P^\gamma_{n,paths}(u)=P^{\psi_\gamma}_n(u^n) \cdot \prod_{k=1}^{n-1} T^{k+1}_k(u^{k+1},u^k), \ \ \ u \in \mb{J}_{n, seq}.
$$
Note that $P^\gamma$ supported on $\mathbb{J}_{\infty,paths} \subset \mb{J}_{\infty, seq}$, and the consistency relation of Proposition \ref{partitionmeasure} guarantee $P^\gamma$ is well-defined. The mapping
$$
\mb{J}_{\infty, seq} \ni u \mapsto  \{ \widetilde{u}^n_k = u^n_k + (r_n-k) \ | \ n \geq 1, \ 1 \leq k \leq r_n \}.
$$ 
pushes $P^\gamma$ forward to determine a probability measure  $\xi^\gamma$ on $2^{\Z_{\geq 0} \times \Z_{> 0}}$, which the last section showed is the fixed-time distribution of the random point configuration $\mc{X}_{(\alpha,\beta)}(\gamma) \subset \Z_{\geq 0} \times \Z_{> 0}$ described in the introduction. We will also refer to this point process as the $(\alpha, \beta)$--\emph{Plancherel point process} (this title does not depend on the choice of coordinates).

\subsection{Computation of correlation kernel} 

For any $n \geq 1$, a point configuration $\mb{X}_n := \{ x^m_k \in \Z | \ 1 \leq m \leq n, \  1 \leq k \leq r_m \}$ in $\Z_{\geq 0} \times [n]$ determines the cylinder set $C_u : = \widetilde{\mc{P}}^{-1}(\mb{X}_n)$ with $u \in \mb{J}_{n,seq}$ satisfying $\widetilde{u}^m_k = x^m_k$, $1 \leq m \leq n$. Using Lemma \ref{indication}, our push-forward measure $\xi^\gamma = P^\gamma \circ \widetilde{\mc{P}}^{-1}$ takes the determinantal form
\begin{equation} \label{detmeasure}
\ba 
\xi^\gamma(\mb{X}_n)=P^\gamma(C_u) = (P_{n,paths}^{\psi_0} A^{\psi_\gamma}_n)(u) =  \prod_{m=1}^{n-1}  \det \left[ \phi_m(x^m_i, x^{m+1}_{j} ) \right]_{i,j=1}^{r_{m+1}} \cdot \det \left [ \Psi^{n, \psi_\gamma}_{r_n-j}(x^{n}_i) \right]_{i,j=1}^{r_n}
\ea
\end{equation}
where $\Psi^{n, \psi_\gamma}_{r_n-i}(s) := \langle \frac{\bar{P}_{s}^{(\alpha_n)}}{2^{\alpha+\beta+1}\Gamma(\alpha+1)}, (x-1)^{r_n-i} R^{\psi_\gamma}_{i-r_n} \rangle_{\alpha_n}$ was defined in Lemma \ref{comprule} (note also $d_1(\lambda) = \bar{\bar{P}}^{(\alpha)}_\lambda(1) = 1$, so there is no constant). Also, define convolution over $\Z_{\geq 0}$ by $(f*g)(x,y): = \sum_{z \geq 0} f(x,z) g(z,y)$ for bivariate functions $f,g$ and by $(f*g)(x):= \sum_{z \geq 0} f(x,z) g(z)$ if $g$ is univariate.
\bp \label{checkconditions}
Let $E(x):=\psi_\gamma(x)$. For any $(s,n), (t,m) \in \Z_{\geq 0} \times \Z_{> 0}$ and $k \in \Z$, define the functions
$$
\ba
 \Phi^{m}_{r_m-k}(t) &:=  \frac{1}{2\pi i} \oint \frac{ \bar{\bar{P}}_{t}^{(\alpha_m)}(w) }{E(w)(w-1)^{r_m-k+1}} dw \\
\phi^{[n,m)}(s,t) &:=-\frac{1}{2 \pi i}  \oint \left \langle \frac{\bar{P}_{s}^{(\alpha_n)}}{2^{\alpha+\beta+1}\Gamma(\alpha+1)}, \frac{\bar{\bar{P}}_{t}^{(\alpha_m)}(u)(u-1)^{r_n-r_m} }{x-u} \right \rangle_{\alpha_n} du, \ \ \  \text{for} \ \  n < m, 
\ea
$$
where the contours are positively oriented (i.e., counterclockwise) simple loops around $[-1,1]$. Then the simple point process $\mc{X}^\gamma_{(\alpha,\beta)}$ determined by $\xi^\gamma$ of \eqref{detmeasure} has determinantal correlation functions $\rho_k^\gamma$ with kernel
\begin{equation} \label{prekernel}
K^\gamma((s,n), (t, m)) = - \phi^{[n,m)}(s,t)1_{(n<m)} + \sum_{k=1}^{r_m} \Psi^{n}_{r_{n} - k}(s) \Phi^{m}_{r_{m} - k}(t).
\end{equation}
\ep
\bpf
The proof relies on linear algebraic details that can be found, e.g., from a combination of Theorem 4.2 of \cite{abpf2} and Lemma 3.4 of \cite{bfs1}, which in turn rely on the Eynard-Mehta theorem in the manner of \cite{br1}. For {\color{black}an exposition} of these details in our setting, {\color{black} the appendix of \cite{cerenzia1} is relevant to this case}, as long as one uses the calculation \eqref{Muppertriangular} to justify that the $r_n \times r_n$ matrix $M$ with ${\color{black} (i,j)}$ entry  
$$
(\phi_{q_i - 1} * \phi_{q_i} * \cdots * \phi_{n-1} * \Psi^{n, \psi_\gamma}_{r_n-j})(-1), \ \ \ q_k := 2k-1,
$$
is upper triangular. Thus, it suffices to identify the functions used there with the functions in the statement above. 

By orthogonal decomposition \eqref{orthogonaldecomposition}, we have for $1 \leq k,l \leq r_n$ 
$$
\ba
\sum_{s \geq 0} \Psi^{n}_{r_n - k}(s)  \Phi^{n}_{r_n-l}(s) &= \frac{1}{2 \pi i} \oint \sum_{s \geq 0} \langle \tilde{P}_{s}^{(\alpha_n)}, (x-1)^{r_n-k} E \rangle_{\alpha_n} \tilde{P}_{s}^{(\alpha_n)}(w) \frac{1}{E(w)(w-1)^{r_n-l+1}} dw \\
& = \frac{1}{2 \pi i} \oint \frac{1}{(w-1)^{k-l+1}} dw = \delta_{kl}.
\ea
$$
Note also that $\Phi^{m}_{r_m-k}(t)$ is a polynomial in $t$ of the same degree as $(\phi_{q_k-1} * \phi^{[q_k,m)})(-1,t)$, where $q_k = 2k-1$. These items confirm that $\{ \Phi^{m}_{r_m-k}(t) \}_{k=1}^{r_m}$ is the unique basis of the linear span of $\{ (\phi_{q_k-1} * \phi^{[q_k,m)})(-1,t) \}_{k=1}^{r_m}$ that is biorthogonal to the $\{\Psi^{n}_{r_n - k}(s)\}_{k=1}^{r_n}$.

We also need to show $\phi^{[n,m)}= \phi_{n} * \cdots * \phi_{m-1}$ for $n<m$. First assume $m = n+1$. If $n$ is odd, then $r_n = r_{n+1}$, $\alpha_n = \alpha$, $\alpha_{n+1} = \alpha+1$, and $\phi_n(s, t) = \phi_n(s) 1_{(s \leq t)}$. Using \eqref{barsidentity}, we have 
$$
\ba
\phi_n(s,t) & = \phi_n(s) 1_{(s \leq t)} = \phi_n(s) \left \langle \tilde{P}^{(\alpha)}_s, \sum_{r=0}^t \tilde{P}^{(\alpha)}_s \right \rangle_{\alpha} =  \frac{1}{2^{\alpha+\beta+1} \Gamma(\alpha+1)} \cdot \left \langle \phi_n(s) \bar{\bar{P}}^{(\alpha)}_s, \sum_{r=0}^t \bar{P}^{(\alpha)}_r \right \rangle_{\alpha} \\
& =  \frac{1}{2^{\alpha+\beta+1} \Gamma(\alpha+1)} \left \langle \bar{P}^{(\alpha)}_s,  \bar{\bar{P}}^{(\alpha+1)}_r \right \rangle_{\alpha} \\
&= -  \frac{1}{2\pi i} \oint  \left \langle \frac{\bar{P}^{(\alpha)}_s}{2^{\alpha+\beta+1} \Gamma(\alpha+1)} ,  \frac{ \bar{\bar{P}}^{(\alpha+1)}_r(u)}{x-u} \right  \rangle_{\alpha} du,
\ea
$$
where the last equality computes the residue at $u = x$, and where we used Proposition \ref{identities} and orthogonality. Similarly, if $n$ is even, $\phi_n(s,t) = \phi_n(s) \cdot 1_{(s<t)}$, $r_n = r_{n+1}-1$, and the $\alpha_\cdot$'s switch, so that
$$
\ba
\phi_n(s,t) & = \phi_n(s) 1_{(s < t)} = \phi_n(s) \left \langle \tilde{P}^{(\alpha+1)}_s, \sum_{r=0}^{t-1} \tilde{P}^{(\alpha+1)}_s \right \rangle_{\alpha+1} =  \frac{1}{2^{\alpha+\beta+1} \Gamma(\alpha+1)} \cdot \left \langle \phi_n(s) \bar{\bar{P}}^{(\alpha+1)}_s, \sum_{r=0}^{t-1} \bar{P}^{(\alpha+1)}_r \right \rangle_{\alpha+1} \\
& =  \frac{1}{2^{\alpha+\beta+1} \Gamma(\alpha+1)} \left \langle \bar{P}^{(\alpha+1)}_s,  \frac{\bar{\bar{P}}^{(\alpha)}_r-1}{x-1} \right \rangle_{\alpha+1}  \\
& =- \frac{1}{2\pi i} \oint  \left \langle \frac{\bar{P}^{(\alpha+1)}_s}{2^{\alpha+\beta+1} \Gamma(\alpha+1)} ,  \frac{ \bar{\bar{P}}^{(\alpha)}_r(u)}{x-u} \right  \rangle_{\alpha+1} \frac{1}{u-1} du
\ea
$$
where now the last equality computes residues at \emph{both} $u = x$ and $u = 1$, and uses the fact that $\bar{\bar{P}}^{(\alpha)}_t(1) = 1$. The full statement for general $n<m$ then follows by induction.
\epf

From our current \eqref{prekernel}, a few more calculations are required to arrive at our main expression \eqref{explicitkernel} for $K^\gamma$ in Theorem \ref{kernel}. 
\bl For any $n,m \geq 1$,
\begin{equation} \label{nearlykernel}
\ba
\sum_{k=1}^{r_m} & \Psi^{n}_{r_{n} - k}(s) \Phi^{m}_{r_{m} - k}(t)
 =  \left \langle \frac{\bar{P}_{s}^{(\alpha_n)}}{2^{\alpha+\beta+1}\Gamma(\alpha+1)}, \frac{ (x-1)^{r_n-r_m}R^E_{r_m-r_n}}{E} \bar{\bar{P}}_{t}^{(\alpha_m)} \right \rangle_{\alpha_n} \\
 & + \frac{1}{2 \pi i} \oint   \left \langle \frac{\bar{P}_{s}^{(\alpha_n)}}{2^{\alpha+\beta+1}\Gamma(\alpha+1)}, \frac{(u-1)^{r_n}(R^E_{r_m-r_n}(u) - E(u)) +(x-1)^{r_n}E}{(x-u)E(u)(u-1)^{r_m}} \right \rangle_{\alpha_n} \bar{\bar{P}}_{t}^{(\alpha_m)}(u) du.\\
\ea
\end{equation}
\el
Note how the terms in \eqref{nearlykernel} simplify if $r_n \geq r_m$, since in this case $R^E_{r_m-r_n} \equiv E$. 

\bpf
If $r_m \leq r_n$, then using the geometric sum identity $\sum_{k=1}^q \left ( \frac{u-1}{x-1} \right )^k = \left ( \frac{u-1}{x-u} \right ) \left ( 1- \left ( \frac{u-1}{x-1}\right )^q \right)$, we have
\begin{equation} \label{geometric}
\ba
\sum_{k=1}^{r_m} \Psi^{n}_{r_{n} - k}(s) \Phi^{m}_{r_{m} - k}(t) = \frac{1}{2 \pi i} \oint \frac{\bar{\bar{P}}_{t}^{(\alpha_m)}}{E(u)} \left \langle \frac{\bar{P}_{s}^{(\alpha_n)}}{2^{\alpha+\beta+1}\Gamma(\alpha+1)}, \frac{(x-1)^{r_n}}{(u-1)^{r_m}} \frac{E(x)}{x-u} \left ( 1- \left ( \frac{u-1}{x-1} \right )^{r_m}  \right ) \right \rangle_{\alpha_n} du.
\ea
\end{equation}
and taking the residue at $u=x$ of the second summand in the inner product yields \eqref{nearlykernel}. If instead $r_m>r_n$, write
$$
\sum_{k=1}^{r_m} \Psi^{n}_{r_{n} - k}(s) \Phi^{m}_{r_{m} - k}(t) = \sum_{k=1}^{r_n} \Psi^{n}_{r_{n} - k}(s) \Phi^{m}_{r_{m} - k}(t) + \sum_{k=r_n+1}^{r_m} \Psi^{n}_{r_{n} - k}(s) \Phi^{m}_{r_{m} - k}(t).
$$
The first summand is treated as in \eqref{geometric}:
$$
\ba
\sum_{k=1}^{r_n} \Psi^{n}_{r_{n} - k}(s) \Phi^{m}_{r_{m} - k}(t) = \frac{1}{2 \pi i} \oint \frac{\bar{\bar{P}}_{t}^{(\alpha_m)}}{E(u)} \left \langle \frac{\bar{P}_{s}^{(\alpha_n)}}{2^{\alpha+\beta+1}\Gamma(\alpha+1)}, \frac{(x-1)^{r_n}}{(u-1)^{r_m}} \frac{E(x)}{x-u} \left ( 1- \left ( \frac{u-1}{x-1} \right )^{r_n}  \right ) \right \rangle_{\alpha_n} du.
\ea
$$
For the second summand, we use the geometric sum identity repeatedly to get
$$
\ba
 \sum_{k=r_n+1}^{r_m}  (E(x) -E(1)) & \frac{(u-1)^k}{(x-1)^{k-r_n}}  - 1_{(r_m \geq r_n+1)}  \sum_{k=r_n+2}^{r_m} \sum_{r=1}^{k-r_n-1} \frac{E^{(r)}(1)}{r!} \frac{(u-1)^k}{(x-1)^{k-r_n-r}}\\
 =  (E(x) &-E(1)) \frac{(u-1)^{r_n+1}}{x-u}  \left ( 1- \left ( \frac{u-1}{x-1}\right )^{r_m-r_n} \right) \\
  &- 1_{(r_m \geq r_n+1)} \sum_{r=1}^{r_m-r_n-1} \frac{E^{(r)}(1)}{r!} \frac{(u-1)^{r_n+r+1}}{x-u} \left ( 1- \left ( \frac{u-1}{x-1}\right )^{r_m-r_n-r} \right).
\ea
$$
and then arrive at
$$
\ba
\sum_{k=r_n+1}^{r_m} & \Psi^{n}_{r_{n} - k}(s) \Phi^{m}_{r_{m} - k}(t) = - \frac{1}{2 \pi i } \oint \bar{\bar{P}}_{t}^{(\alpha_m)}(u) \left \langle \frac{\bar{P}_{s}^{(\alpha_n)}}{2^{\alpha+\beta+1}\Gamma(\alpha+1)}, \frac{(x-1)^{r_n-r_m} R^E_{r_m-r_n}}{E(u)(x-u)} \right \rangle_{\alpha_n} du  \\
& + \frac{1}{2 \pi i } \oint \bar{\bar{P}}_{t}^{(\alpha_m)}(u) \left \langle \frac{\bar{P}_{s}^{(\alpha_n)}}{2^{\alpha+\beta+1}\Gamma(\alpha+1)}, \frac{(u-1)^{r_n-r_m} (R^E_{r_m-r_n}(u)-E(u)+E(x))}{E(u)(x-u)} \right \rangle_{\alpha_n} du
\ea
$$
Taking the residue at $u=x$ of the first term and combining with the first summand completes the proof.
\epf
Now plug the expression \eqref{nearlykernel} into \eqref{prekernel}. If $r_n \geq r_m$, the expression \eqref{nearlykernel} simplifies, quickly leading to the main expression \eqref{explicitkernel} for $K^\gamma$ in Theorem \ref{kernel} (note we have multiplied by the conjugating factor ``$(-1)^{r_n-r_m}$", which vanishes in the determinant). But if $r_n<r_m$, we get \eqref{explicitkernel} along with the additional term
$$
\ba
\frac{1}{2 \pi i} \oint & \left \langle \frac{\bar{P}_{s}^{(\alpha_n)}}{2^{\alpha+\beta+1}\Gamma(\alpha+1)}, \bar{\bar{P}}_{t}^{(\alpha_m)}(u) \frac{(u-1)^{r_n-r_m}R^E_{r_m-r_n}(u)}{E(u)(x-u)} \right \rangle_{\alpha_n} du \\
& + \left \langle \frac{\bar{P}_{s}^{(\alpha_n)}}{2^{\alpha+\beta+1}\Gamma(\alpha+1)}, \bar{\bar{P}}_{t}^{(\alpha_m)} \frac{(x-1)^{r_n-r_m}R^E_{r_m-r_n}}{E} \right \rangle_{\alpha_n}
\ea
$$
The residue of the first term at $u=x$ exactly cancels the second, and since $$(u-1)^{r_n-r_m}R^E_{r_m-r_n}(u) = \sum_{k=r_m-r_n}^\infty \frac{E^{(k)}(1)}{k!}(u-1)^{k+r_n-r_m},$$ there is no residue at $u=1$ (note of course $E=\psi_\gamma$ does not have a pole at $1$ either). This completes the proof of Theorem \ref{kernel}. 

\section{Edge asymptotic analysis} \label{sectionasymptotics}

\subsection{Finite distance from the wall: Jacobi kernel}

Assume again that $\gamma \sim N \cdot \tau >0$ and $r_{n}, r_m \sim N \cdot \eta>0 $, but that $s, t$ are fixed and finite. Assume only the difference $n - m$ is of constant order. Let $A(z) := \tau z + \eta \log (1-z)$ and note $1-\eta/\tau$ is a zero of $A'(z)$. Write our kernel as 
$$
\ba
  \frac{1}{2 \pi i}& \int_{-1}^1 \oint \frac{e^{- N(A(1-\eta/\tau) - A(x))}}{e^{ N(A(u) - A(1-\eta/\tau))}}  \frac{\bar{P}_{s}^{(\alpha_n)}(x)}{2^{\alpha+\beta+1}\Gamma(\alpha+1)} \bar{\bar{P}}_{t}^{(\alpha_m)}(u) \frac{(1-x)^{\alpha_n}(1+x)^{\beta}}{x-u} dudx\\
& +1_{(n\geq m)} \left \langle \frac{\bar{P}_{s}^{(\alpha_n)}}{2^{\alpha+\beta+1}\Gamma(\alpha+1)}, (x-1)^{r_n-r_m} \bar{\bar{P}}_{t}^{(\alpha_m)} \right \rangle_{\alpha_n}. 
\ea
$$
\begin{figure}[h]
\centering
\includegraphics[width=0.45\textwidth]{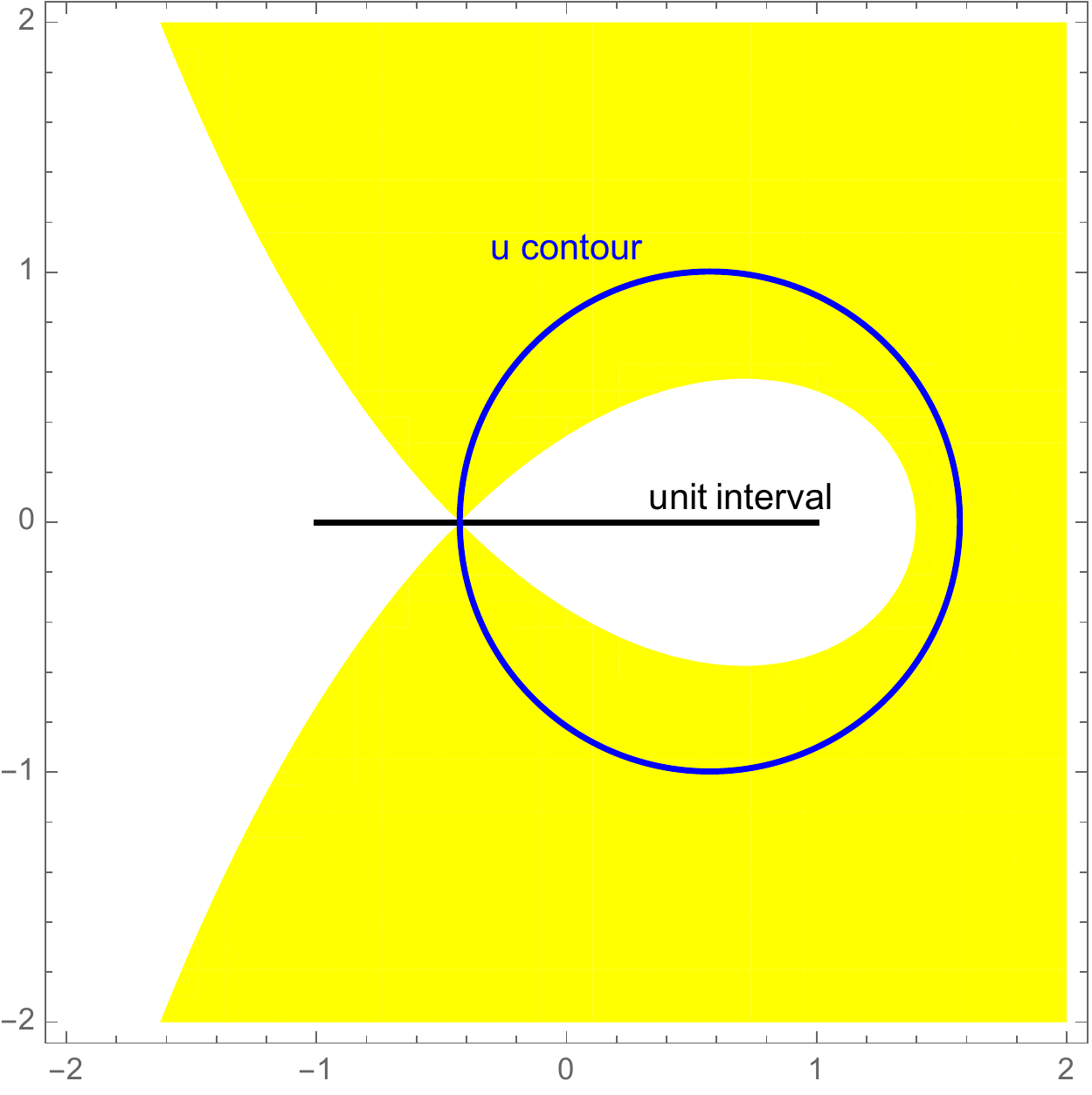}
\caption{The plot exemplifies steepest descent deformations when targeting a finite distance from the wall. The yellow region signifies where $\mf{R}(A(u)-A(1-\eta/\tau))>0$ and the white elsewhere.} \label{jacobiasymptoticsFIGURE}
\end{figure}
Now deform the $u$-contour, as in Figure \ref{jacobiasymptoticsFIGURE}, to be a steepest ascent loop remaining in the region $ \mc{R}(A(u) - A(1-\eta/\tau)) > 0$ and passing through the critical point $1-\eta/\tau$. If $1-\eta/\tau< -1$, the unit interval already lies in the region $ \mc{R}(A(1-\eta/\tau) - A(x)) >0$, so the double integral term tends to zero without picking up residues. As for the frozen region above, the kernel converges to a triangular matrix whose diagonal entries are $1$'s. But if $1-\eta/\tau> -1$, then we acquire a residue at $u = x$ given by
$$
 -  \frac{(-1)^{r_n-r_m}}{2^{\alpha+\beta+1}\Gamma(\alpha+1)}  \int_{-1}^{1-\eta/\tau} \bar{P}_{s}^{(\alpha_n)}(x) \bar{\bar{P}}_{t}^{(\alpha_m)}(x) (1-x)^{r_n-r_m+\alpha_n} (1+x)^{\beta} dx
$$
Noting that $(-1)^{r_n-r_m}$ is a conjugating factor completes the proof of Theorem \ref{jacobilimit}.

\subsection{Hard--edge Pearcey kernel}

Assume that $\gamma \sim N/2$, and $s \sim N^{1/4} \nu_1>0$, $t \sim N^{1/4} \nu_2>0$, and lastly that $r_{n} - N\sim \sqrt{N} \sigma_1$, $r_{m} - N\sim \sqrt{N} \sigma_2$. Assume also that $n,m$ have constant parity. The Mehler--Heine formula (Chapter 8 of Szeg{\color{black} \H{o}} \cite{gs1}) {\color{black} and the symmetry relation $P_n^{(\alpha,\beta)}(-z) = (-1)^n P_n^{(\beta,\alpha)}(z)$} tells us that for any $\alpha, \beta > -1$,
\begin{equation} \label{MHformula}
\lim_{N \to \infty} N^{-\beta/4} (-1)^{\lfloor \nu_1 N^{1/4} \rfloor } P^{{\color{black} (\alpha, \beta)}}_{\lfloor \nu_1 N^{1/4} \rfloor} \left ( \frac{z}{\sqrt{N}}-1\right)  =  \left ( \frac{\sqrt{2z}}{2} \right )^{-\beta} J_\beta(\nu_1\sqrt{2z}).
\end{equation}
{\color{black} Note that the right--hand--side only depends on $\beta$.} Recall the {\color{black} particle--hole involution} $(\Z_{ \geq 0} \times \Z_{\geq 1}) \setminus (\mc{X}(\gamma))$ of $\mc{X}^\gamma_{(\alpha,\beta)}$ is also determinantal with correlation kernel 
\begin{equation} \label{compkernel}
\ba
K^\gamma_\Delta((s,n),(t,m)) &:= \delta_{s,t} \delta_{n,m} - K^\gamma((s,n),(t,m)) \\
 &=-  \frac{1}{2 \pi i}  \int_{-1}^1 \oint \frac{\psi_\gamma(x)}{\psi_\gamma(u)} \frac{\bar{P}_{s}^{(\alpha_n)}(x)}{2^{\alpha+\beta+1}\Gamma(\alpha+1)} \bar{\bar{P}}_{t}^{(\alpha_m)}(u) \frac{(x-1)^{r_n}}{(u-1)^{r_m}} \frac{(1-x)^{\alpha_n}(1+x)^{\beta}}{x-u} dudx\\
& -1_{(n>m)} \left \langle \frac{\bar{P}_{s}^{(\alpha_n)}}{2^{\alpha+\beta+1}\Gamma(\alpha+1)}, (x-1)^{r_n-r_m} \bar{\bar{P}}_{t}^{(\alpha_m)} \right \rangle_{\alpha_n}.
\ea
\end{equation}
where we have used the fact that
$$
1_{(n=m)} \delta_{s,t} = 1_{(n=m)} \left \langle \frac{\bar{P}_{s}^{(\alpha_n)}}{2^{\alpha+\beta+1}\Gamma(\alpha+1)}, (x-1)^{r_n-r_m} \bar{\bar{P}}_{t}^{(\alpha_m)} \right \rangle_{\alpha_n}
$$
\begin{figure}[h]
\centering
\includegraphics[width=0.45\textwidth]{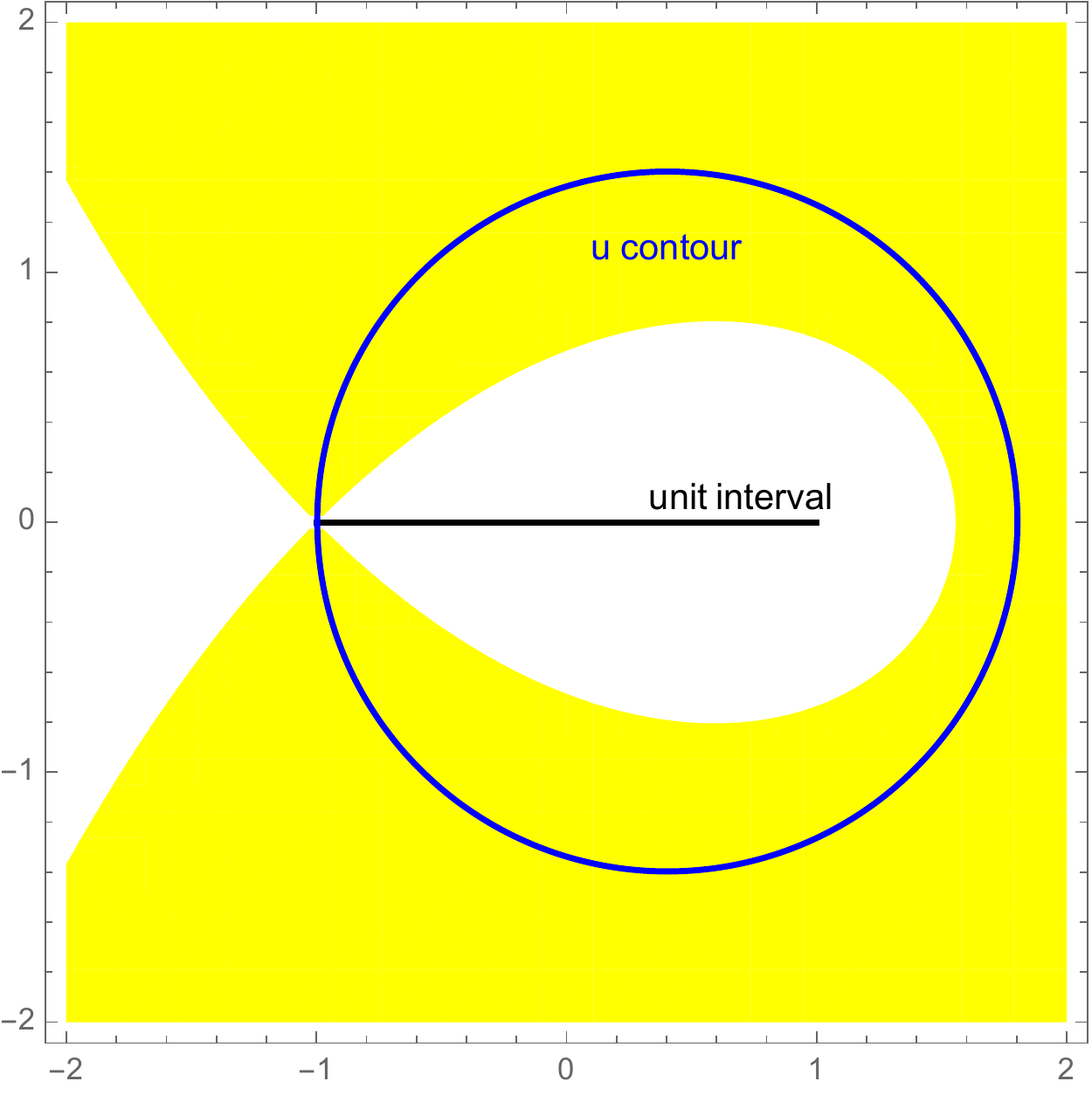}
\caption{Steepest descent deformations. The yellow region signifies the region $\mf{R}(A(u)-A(-1))>0$ and the white elsewhere.} \label{pearceyasymptotics}
\end{figure}
Deforming the $u$-contour in the double integral term of \eqref{compkernel} as in Figure \ref{pearceyasymptotics}, we endeavor to find the contribution at $-1$. Making the substitutions $x' = N^{1/2}(x+1)$ and $u' = N^{1/2}(u+1)$, we have, for large $N$,
$$
\ba
\frac{(1-x)^{\alpha_n}(1+x)^{\beta}}{x-u} dudx & = N^{1/2} \frac{(2-x'N^{-1/2})^{\alpha_n}(x')^{\beta} N^{-\beta/2}}{x'-u'} \cdot N^{-1} du'dx' \\
& \sim N^{-1/2-\beta/2} \frac{2^{\alpha_n}(x')^\beta}{x'-u'} \cdot du'dx'.
\ea
$$
Putting this together with the Mehler--Heine formula \eqref{MHformula}, we get
$$
\ba
&(-1)^{s-t}  \cdot \frac{\bar{c}_s^n \bar{\bar{c}}_t^m}{2^{\alpha+\beta+1}\Gamma(\alpha+1)} P_{s}^{(\alpha_n)}(x) P_{t}^{(\alpha_m)}(u)\frac{(1-x)^{\alpha_n}(1+x)^{\beta}}{x-u} dudx \\
& \sim N^{-1/2} \left ( \frac{x'}{u'} \right )^{\beta/2} \frac{J_\beta(\nu_1\sqrt{2x'}) J_\beta(\nu_2\sqrt{2u'}) }{x'-u'} du'dx'
\begin{cases}
(\nu_2 N^{1/4})^{\alpha+1} (\nu_1 N^{1/4})^{-\alpha}, & n \  \text{even}, \ m \  \text{even} \\ 
\Gamma(\alpha+1) (\nu_2 N^{1/4})^{-\alpha} (\nu_1 N^{1/4})^{-\alpha}, & n \  \text{even}, \ m \  \text{odd} \\ 
\frac{1}{\Gamma(\alpha+1)} (\nu_2 N^{1/4})^{\alpha+1} (\nu_1 N^{1/4})^{\alpha+1}, & n \  \text{odd}, \ m \  \text{even} \\ 
(\nu_2 N^{1/4})^{-\alpha} (\nu_1 N^{1/4})^{\alpha+1} , & n \  \text{odd}, \ m \  \text{odd} \\ 
\end{cases} \\
& =\frac{(N^{1/4})^{(2 \alpha+1)1_{(n \ odd)} - (2\alpha+1)1_{(m \ odd)}}}{\Gamma(\alpha+1)^{1_{(m \ \text{even})} -1_{(n \ \text{even})}}} \frac{(\nu_1)^{(-1)^{n+1}(\alpha+1/2)}}{(\nu_2)^{(-1)^{m+1}(\alpha+1/2)}} \cdot N^{-1/4} \sqrt{\nu_1 \nu_2} \left ( \frac{x'}{u'} \right )^{\beta/2} \frac{J_\beta(\nu_1\sqrt{2x'}) J_\beta(\nu_2\sqrt{2u'}) }{x'-u'} du'dx'
\ea
$$
where we have omitted the rest of the series ``$s^{a-b} \sum_{k=1}^\infty a_k s^{-k}$" because the conjugating factor $(N^{1/4})^{(2 \alpha+1)1_{(n \ odd)} - (2\alpha+1)1_{(m \ odd)}}$ drops out when taking the determinant. The remainder of the integrand satisfies
$$
\ba
(-1)^{r_m-r_n} 2^{r_m-r_n} \frac{e^{\gamma x}}{e^{\gamma u}}  \frac{(x-1)^{r_n}}{(u-1)^{r_m}} & = \frac{ e^{-N(A(-1)-A(x))+ \sigma_1 \sqrt{N}(\log(1-x) - \log 2)}}{e^{N(A(u) - A(-1))+ \sigma_2 \sqrt{N}(\log(1-u) - \log 2)}} \sim \frac{ e^{-(x')^2/8 - \sigma_1 x'/2}}{ e^{-(u')^2/8 - \sigma_2 u'/2}},
\ea 
$$
where $A(z): = z/2+ \log(1-z)$ and the last approximation follows from, respectively, second and first order Taylor expansions about $-1$. Hence, we have (omitting conjugating factors)
$$
\ba
&  N^{1/4} \frac{1}{2 \pi i}  \int_{-1}^1 \oint \frac{e^{\gamma x}}{e^{\gamma u}} \frac{\bar{P}_{s}^{(\alpha_n)}(x)}{2^{\alpha+\beta+1}\Gamma(\alpha+1)} \bar{\bar{P}}_{t}^{(\alpha_m)}(u) \frac{(x-1)^{r_n}}{(u-1)^{r_m}} \frac{(1-x)^{\alpha_n}(1+x)^{\beta}}{u-x} dudx \\
& \to  \frac{1}{2 \pi i}  {\color{black}\int_0^{\infty} \int_{-i\infty}^{i\infty}  } e^{ \frac{(u')^2-(x')^2}{8} + \frac{\sigma_2 u' - \sigma_1 x'}{2}} \left ( \frac{x'}{u'} \right )^{\beta/2}  \sqrt{\nu_1 \nu_2} J_\beta(\nu_1\sqrt{2x'})  J_\beta(\nu_2\sqrt{2u'}) \frac{du'dx'}{x'-u'}
\ea
$$
Turning to the other term of $K^\gamma_\Delta$, note the following two asymptotic relations:
$$
(-2)^{r_m-r_n} (x-1)^{r_n-r_m} = \left (1-\frac{x'}{2 \sqrt{N}} \right )^{(\sigma_1 - \sigma_2) \sqrt{N}} \to e^{\frac{-(\sigma_1-\sigma_2)x'}{2}}.
$$
$$
\ba
 & \frac{\bar{c}^n_s \bar{\bar{c}}^m_t (-1)^{s-t}}{2^{\alpha+\beta+1}\Gamma(\alpha+1)} P_{s}^{(\alpha_n)}(x) P_{t}^{(\alpha_m)}(x) (1-x)^{\alpha_n}(1+x)^{\beta} dx \\
 &  \sim \frac{(N^{1/4})^{(2 \alpha+1)1_{(n \ odd)} - (2\alpha+1)1_{(m \ odd)}}}{\Gamma(\alpha+1)^{1_{(m \ \text{even})} -1_{(n \ \text{even})}}} \frac{(\nu_1)^{(-1)^{n+1}(\alpha+1/2)}}{(\nu_2)^{(-1)^{m+1}(\alpha+1/2)}} \cdot \sqrt{\nu_1 \nu_2} N^{-1/4} J_\beta(\nu_1\sqrt{2x'}) J_\beta(\nu_2\sqrt{2x'})  dx'
\ea
$$
Then the single-integral term satisfies
$$
\ba
 \frac{(-2)^{r_m-r_n}  \Gamma(\alpha+1)^{1_{(m \ \text{even})} -1_{(n \ \text{even})}}}{(N^{1/4})^{(2 \alpha+1)1_{(n \ odd)} - (2\alpha+1)1_{(m \ odd)}}} & \times N^{1/4}1_{(n>m)} \left \langle \frac{\bar{P}_{s}^{(\alpha_n)}}{2^{\alpha+\beta+1}\Gamma(\alpha+1)}, (x-1)^{r_n-r_m} \bar{\bar{P}}_{t}^{(\alpha_m)} \right \rangle_{\alpha_n}  \\
&  \to 1_{(\sigma_1 > \sigma_2)} \left ( \int_0^\infty  e^{\frac{-(\sigma_1-\sigma_2)x'}{2}} \sqrt{\nu_1 \nu_2} J_\beta(\nu_1\sqrt{2x'}) J_\beta(\nu_2\sqrt{2x'}) dx' \right ) 
\ea
$$

It remains to show that the expression for the kernel here matches the hard--edge Pearcey kernel $L^{\beta}$. By the identity (see e.g. equation (25) of \cite{bmp}),
$$
\int_0^{\infty} J_{\beta}(a_1 x) J_{\beta}(a_2 x) e^{-\gamma^2 x^2}x dx 
= \frac{1}{2\gamma^2} \exp\left( - \frac{ a_1^2 + a_2^2 }{4\gamma^2} \right) I_{\beta}\left( \frac{a_1a_2}{2\gamma^2} \right), \quad \mathrm{Re}\ \beta>-1, \quad \mathrm{Re}\ \gamma^2 >0,
$$
with $x'=2x^2, \gamma^2 = (\sigma_1-\sigma_2), a_1 = 2\nu_1, a_2 = 2\nu_2,$ then
\begin{multline*}
1_{(\sigma_1 > \sigma_2)} \left ( \int_0^\infty  e^{\frac{-(\sigma_1-\sigma_2)x'}{2}} \sqrt{\nu_1 \nu_2} J_\beta(\nu_1\sqrt{2x'}) J_\beta(\nu_2\sqrt{2x'}) dx' \right ) \\
=1_{(\sigma_1 > \sigma_2)} \sqrt{\nu_1\nu_2} \frac{2}{\sigma_1-\sigma_2} \exp\left( -\frac{\nu_1^2+\nu_2^2}{\sigma_1-\sigma_2}\right) I_{\beta} \left( \frac{2\nu_1\nu_2}{\sigma_1-\sigma_2}\right).
\end{multline*}
By the substitutions $x'=2x^2,u'=2u^2,$
\begin{multline*}
\frac{1}{2 \pi i}  {\color{black}\int_0^{\infty} \int_{-i\infty}^{i\infty}  } e^{ \frac{(u')^2-(x')^2}{8} + \frac{\sigma_2 u' - \sigma_1 x'}{2}} \left ( \frac{x'}{u'} \right )^{\beta/2}  \sqrt{\nu_1 \nu_2} J_\beta(\nu_1\sqrt{2x'})  J_\beta(\nu_2\sqrt{2u'}) \frac{du'dx'}{x'-u'} \\
= \sqrt{\nu_1\nu_2} \frac{4}{\pi i} \int_0^{\infty} \int_C  e^{ \tfrac{u^4 - x^4 }{2} + \sigma_2 u^2 - \sigma_1 x^2} \left( \frac{x}{u} \right)^{\beta} J_{\beta}(2\nu_1 x) J_{\beta}(2\nu_2 u)  xu \frac{dudx}{x^2 - u^2},
\end{multline*}
where $C$ consists of two rays, one from $e^{i \pi /4} \infty$ to $0$ and one from $0$ to $e^{-i \pi / 4}\infty$. The final expression is therefore
$$
L^{\beta}(\sigma_1,\nu_1^2;\sigma_2,\nu_2^2) \cdot 2\sqrt{\nu_1\nu_2},
$$
finishing the proof of Theorem \ref{pearceylimit}.

{\color{black} 
 
\begin{remark}\label{CoV} The expression $2\sqrt{\nu_i\nu_j}$ can be understood with the following heuristics. If $\mathfrak{X}$ is a determinantal point process on $\mathbb{R}$ with correlation kernel $K(\cdot,\cdot)$ and correlation functions $\rho_k$, then (see e.g. (2.6) and (2.16) of \cite{kj1})
\begin{align*}
\mathbb{E}[ \# \text{ of pairs of distinct points of } \mathfrak{X} \text{ in } A \subset \mathbb{R}] &= \int_{A\times A} \rho_2(x,y)dxdy \\
&= \int_{A\times A} \left( K(x,x)K(y,y)-K(x,y)K(y,x)\right)dxdy.
\end{align*}
If $\phi:\mathbb{R}\rightarrow \mathbb{R}$ is a differentiable bijection, and $\phi(\mathfrak{X})$ is a determinantal point process with correlation kernel $\tilde{K}(\cdot,\cdot)$, then
\begin{multline*}
\mathbb{E}[ \# \text{ of pairs of distinct points of } \tilde{\mathfrak{X}} \text{ in } \phi(A)] = \int_{\phi(A \times A)} \left(\tilde{K}(\tilde{x},\tilde{x})\tilde{K}(\tilde{y},\tilde{y}) -  \tilde{K}(\tilde{x},\tilde{y})\tilde{K}(\tilde{y},\tilde{x}) \right) d\tilde{x} d\tilde{y} \\
= \int_{A\times A} \left( \tilde{K}(\phi(x),\phi(x)) \tilde{K}(\phi(y),\phi(y)) - \tilde{K}(\phi(x),\phi(y)) \tilde{K}(\phi(y),\phi(x)) \right) \phi'(x) \phi'(y) dx dy
\end{multline*}
where the last equality is due to the change of variables $\tilde{x}=\phi(x),\tilde{y}=\phi(y)$. These two integrals must be equal, so $K$ and $\tilde{K}$ are related by
$$
K(x,y) = \tilde{K}(\phi(x),\phi(y)) \sqrt{ \phi'(x) \phi'(y)}.
$$
Note that this can be written equivalently as
$$
\sqrt{ \frac{\phi'(y)}{\phi'(x)} } K(x,y) = \tilde{K}(\phi(x),\phi(y)) \phi'(y),
$$
which is how it is expressed in Corollary 2.25 of \cite{sdbv1}.
\end{remark}
}

\bibliographystyle{plain}
\bibliography{Jacobigrowth}

\end{document}